\documentclass[12pt]{article}

\usepackage{amsmath,amssymb,bm,amsthm,mathrsfs,amscd}
\usepackage{color}
\usepackage{graphicx}%You have to use it when you upload this file to arXiv and comment out \usepackage[dvipdfmx]{graphicx}.
\theoremstyle{definition}
\newtheorem{Def}{Def}[section]

\newtheorem{Them}[Def]{Theorem}
\newtheorem{Lem}[Def]{Lemma}
\newtheorem{Cor}[Def]{Corollary}
\newtheorem{Prop}[Def]{Proposition}
\newtheorem{Examp}[Def]{Example}

\numberwithin{equation}{section}

\title{High order free hyperplane arrangements in $3$-dimensional vector spaces}
\author{Norihiro Nakashima\footnote{Department of Mathematics, Nagoya Institute of Technology, Aichi, 466-8555, Japan. Email: nakashima@nitech.ac.jp }}

\date{}
\begin{document}

\maketitle

\begin{abstract}
Holm introduced $m$-free $\ell$-arrangements which is a generalization of free arrangements, while he asked whether all $\ell$-arrangements are $m$-free for $m$ large enough.
Recently Abe and the author gave a negative answer to this question when $\ell\geq 4$.
In this paper we verify that $3$-arrangements $\mathscr{A}$ are $m$-free and compute the $m$-exponents for all $m\geq |\mathscr{A}|-2$, where $|\mathscr{A}|$ is the cardinality of $\mathscr{A}$.
Hence Holm's question has a positive answer when $\ell=3$.
Finally we prove that $3$-dimensional Weyl arrangements of types A and B are $m$-free for all $m\geq 0$.

\vspace{2mm}

\end{abstract}

\section{Introduction}
Let $\mathbb{K}$ be a field of characteristic zero and
let $V$ be an $\ell$-dimensional vector space over $\mathbb{K}$.
A {\bf (central hyperplane) arrangement} $\mathscr{A}$,
or $(\mathscr{A},V)$, is a finite set of hyperplanes
in $V$ which contain the origin.
We call $\mathscr{A}$ an $\ell$-arrangement
when we emphasize the dimension of $V$.

For a $k$-dimensional vector space $\Omega$ and for a $\mathbb{K}$-basis $\{\omega_1,\dots,\omega_k\}$ for $\Omega$, a $\mathbb{K}$-basis $\{\omega_{1}^{\ast},\dots,\omega_{k}^{\ast}\}$ for $\Omega^{\ast}:={\rm Hom}_{\mathbb{K}}(\Omega,\mathbb{K})$ is said to be the dual basis for $\{\omega_1,\dots,\omega_k\}$ if $\omega_{i}^{\ast}(\omega_{j})=\left\{
\begin{array}{ll}
1&\ (i=j),\\
0&\ (i\neq j).
\end{array}
\right.$
Let $S={\rm Sym}(V^{\ast})$ be the symmetric algebra of $V^{\ast}$,
let $\{x_{1},\dots,x_{\ell}\}$ be a $\mathbb{K}$-basis for $V^{\ast}$, and
let $\{\partial_{1},\dots,\partial_{\ell}\}$ be the dual basis
for $\{x_{1},\dots,x_{\ell}\}$.
We can consider $S$ as the polynomial ring
$\mathbb{K}[x_{1},\dots,x_{\ell}]$ in variables $x_{1},\dots,x_{\ell}$
and $\partial_{1}=\partial/\partial x_{1},\dots,\partial_{\ell}
=\partial/\partial x_{\ell}$ as partial derivatives.
For any hyperplane $H\in\mathscr{A}$,
there exists a linear form $\alpha_{H}$ in the dual space $V^{\ast}$
such that $\{\alpha_{H}=0\}=H$.
We call $Q:=Q(\mathscr{A}):=\prod_{H\in\mathscr{A}}\alpha_{H}$
a defining polynomial of $\mathscr{A}$.
The cardinality $n:=|\mathscr{A}|$ of the arrangement $\mathscr{A}$
equals the degree of $Q$.

Let $\mathbb{N}=\{0,1,2,\dots\}$ be the set of nonnegative integers
and let $m\in\mathbb{N}$.
For $\bm{a}=(a_1,\dots,a_{\ell})\in\mathbb{N}^{\ell}$, we write
\begin{align}
|\bm{a}|=a_1 +\cdots+a_{\ell},\ 
\bm{a}!=a_1 !\cdots a_{\ell}!,\ 
x^{\bm{a}}=x_1^{a_1}\cdots x_{\ell}^{a_{\ell}},\ {\rm and}\ 
\partial^{\bm{a}}=\partial_{1}^{a_{1}}\cdots\partial_{\ell}^{a_{\ell}}.
\end{align}
Let $D^{(m)}(S):=\sum_{|\bm{a}|=m}S\partial^{\bm{a}}$
be the $S$-module generated by $m$-th partial derivatives.
An $S$-submodule $D^{(m)}(\mathscr{A},V)$ of $D^{(m)}(S)$ is defined by
\begin{align}
D^{(m)}(\mathscr{A},V):=\left\{\theta\in D^{(m)}(S)\mid
\theta(QS)\subseteq QS\right\},
\end{align}
which is called the {\bf module of $m$-th order
$\mathscr{A}$-differential operators}. Let
\begin{align}
s_m(\ell):=\binom{m+\ell-1}{m}=\binom{m+\ell-1}{\ell-1}.
\end{align}
We say that an arrangement $\mathscr{A}$ is {\bf $m$-free}
if $D^{(m)}(\mathscr{A},V)$ has linearly independent generators
$\theta_1,\dots,\theta_{s_m(\ell)}$ over $S$, and
the set $\{\theta_1,\dots,\theta_{s_m(\ell)}\}$
of such generators is called a free basis.
Let $S_i$ be the vector space consisting
of homogeneous polynomials of degree $i$ in $S$.
For $\theta=\sum_{|\bm{a}|=m}f_{\bm{a}}\partial^{\bm{a}}\in D^{(m)}(S)$,
we write $\deg(\theta)=i$ if
$f_{\bm{a}}\in S_i$ for each $\bm{a}$.
A multi-set ${\rm exp}_m(\mathscr{A},V)$ of {\bf $m$-exponents}
of an $m$-free arrangement $\mathscr{A}$ is defined by
\begin{align}
{\rm exp}_m(\mathscr{A},V):=\{\deg(\theta_1),\dots,\deg(\theta_{s_m(\ell)})\}
=\{e_1^{i_1},e_2^{i_2},\dots\},
\end{align}
where $e^i\in{\rm exp}_m(\mathscr{A},V)$
means that the integer $e$ occurs $i$ times
in the multi-set ${\rm exp}_m(\mathscr{A},V)$.
We usually call $1$-free arrangements free arrangements.
Free arrangements have been studied, relating with combinatorics
of hyperplane arrangements.
In particular the addition-deletion theorems \cite{Terao-add-del}
and Terao's factorization theorem \cite{Terao-fac}
describe relations of a free arrangement and the intersection poset
\begin{align}
L(\mathscr{A}):=\left\{\bigcap_{H\in\mathscr{B}}H\,\middle|\,\mathscr{B}\subseteq\mathscr{A}\right\}
\end{align}
whose ordering is defined by reverse inclusions (see also \cite{Orlik-Terao}).

For $m\geq 2$, the $m$-free arrangements are introduced by Holm \cite{Holm-phd,Holm} to study the ring of differential operators on $\mathscr{A}$, which is a quotient ring of the Weyl algebra (see \cite{Mac-Rob}).
In particular Holm \cite{Holm-phd} proved that $2$-arrangements are $m$-free for all $m\geq 0$ by constructing free bases, and these free bases are used in \cite{Nakashima-noeth} to prove that the ring of differential operators on a $2$-arrangement is a Noetherian ring.
It is also an interesting problem to observe the behavior of $m$-freeness.
Holm \cite{Holm-phd} asked whether all arrangements are $m$-free for $m$ large enough.
In the case when $\mathscr{A}$ is a generic $\ell$-arrangement (i.e., every $\ell$ hyperplanes of $\mathscr{A}$ intersect only at the origin), $\mathscr{A}$ is $m$-free if and only if $m\geq n-\ell+1$ \cite{Holm-phd,NOS}.
This means that Holm's question has a positive answer for generic arrangements.
However for any $\ell\geq 4$, there exist an $\ell$-arrangement $\mathscr{A}$ such that $\mathscr{A}$ is not $m$-free for any $m\geq 1$ (shown by Abe and the author \cite{Abe-Nakashima}), that is, Holm's question has a negative answer for $\ell\geq 4$.
In this paper we verify that Holm's question has a positive answer for $3$-arrangements by proving that all $3$-arrangements are $m$-free for all $m\geq n-2$.
After that we compute the $m$-exponents of $3$-arrangements for all $m\geq n-2$ which depend only on the structure of the poset $L(\mathscr{A})$\footnote{In other words, let $(\mathscr{A},V)$ and $(\mathscr{A}^{\prime},V)$ be $3$-arrangements with $n=|\mathscr{A}|=|\mathscr{A}^{\prime}|$, and let $m\geq n-2$. If $L(\mathscr{A})\simeq L(\mathscr{A}^{\prime})$ as posets, then ${\rm exp}_m(\mathscr{A},V)={\rm exp}_m(\mathscr{A}^{\prime},V)$.}.
In addition we verify that $3$-dimensional Weyl arrangements of types A and B are $m$-free for all $m\geq 0$.
Since $3$-dimensional irreducible Weyl arrangements are only of types A and B, all the Weyl $3$-arrangements are $m$-free for all $m\geq 0$.

The paper is structured as follows.
In Section \ref{sec-change-vari} we argue a change of variables of $S$
and we define submodules of $D^{(m)}(\mathscr{A},V)$.
In Section \ref{sec-when-m>n-2} we prove that if $m\geq n-2$, then $D^{(m)}(\mathscr{A},V)$ is a direct sum of deformed submodules defined in Section \ref{sec-change-vari}.
By using this direct sum, we prove that $3$-arrangements are $m$-free for all $m\geq n-2$.
In Section \ref{sec-exp-3arr} we compute the $m$-exponents of $3$-arrangements for all $m\geq n-2$.
Finally, in Section \ref{sec-wyle-arrAB}, we construct free bases for $D^{(m)}(\mathscr{A},V)$ for the remaining orders $m<n-2$ when $\mathscr{A}$ is the Weyl $3$-arrangement of type A and of type B.

Data sharing not applicable to this article as no datasets were generated or analysed during the current study.

\section{Change of variables}\label{sec-change-vari}
For any arrangement $\mathscr{A}$, let
\begin{align}
L(\mathscr{A})_1:=
\left\{X\in L(\mathscr{A})\,\middle|\,\dim_{\mathbb{K}}(X)=1\right\}
\end{align}
be the set of $1$-dimensional elements in $L(\mathscr{A})$.
Let $X\in L(\mathscr{A})_1$ and let $v_X$ be a nonzero vector in $X$.
Then $X=\mathbb{K}v_X$.
We define a derivation $\delta_X\in\sum_{i=1}^{\ell}\mathbb{K}\partial_i$ by
\begin{align}
\delta_X:=\sum_{i=1}^{\ell}x_i(v_X)\partial_i.
\end{align}
For $X\in L(\mathscr{A})$, the localization $\mathscr{A}_X$
of $\mathscr{A}$ at $X$ is defined by
\begin{align}
\mathscr{A}_X:=\left\{H\in\mathscr{A}\,\middle|\,X\subseteq H\right\}.
\end{align}
\begin{Prop}\label{prop-deltaX}
For $H\in\mathscr{A}$,
$\delta_X(\alpha_H)=0$
if and only if $H\in\mathscr{A}_X$.
Moreover
\begin{align}\label{eq-deltaX}
\mathbb{K}\delta_X=
\left\{\delta\in\sum_{i=1}^{\ell}\mathbb{K}\partial_i\,\middle|\,
\delta\left(\alpha_H\right)=0\ {\rm for}\ H\in\mathscr{A}_X
\right\}.
\end{align}
\end{Prop}
\noindent
{\it Proof.}
For $H\in\mathscr{A}$, we have
$\alpha_H=\sum_{i=1}^{\ell}\partial_i(\alpha_H)x_i$ and
\begin{align*}
\delta_X(\alpha_H)=\sum_{i=1}^{\ell}x_i(v_X)\partial_i(\alpha_H)
=\sum_{i=1}^{\ell}\partial_i(\alpha_H)x_i(v_X)=\alpha_H(v_X).
\end{align*}
Therefore
$H\in\mathscr{A}_X\ \Leftrightarrow\ 
X\subseteq H\ \Leftrightarrow\ \alpha_H(v_X)=0
\ \Leftrightarrow\ \delta_X(\alpha_H)=0$.
This implies that $\delta_X$ is in the right hand side of \eqref{eq-deltaX}.

It remains to prove that the dimension of
the right hand side of \eqref{eq-deltaX} equals one.
Let $\{e_1,\dots,e_{\ell}\}$ be a basis for $V$
for which $\{x_1,\dots,x_{\ell}\}$ is the dual basis.
Since $\sum_{i=1}^{\ell}\alpha_H(e_i)x_i=\alpha_H=
\sum_{i=1}^{\ell}\partial_i(\alpha_H)x_i$,
we have $\alpha_H(e_i)=\partial_i(\alpha_H)$ for $1\leq i\leq \ell$.
Hence the assertion follows from
\begin{align*}
&\left\{(v_1,\dots,v_{\ell})\in\mathbb{K}^{\ell}\,\middle|\,
\sum_{i=1}^{\ell}v_i\partial_i
\left(\alpha_H\right)=0\ {\rm for}\ H\in\mathscr{A}_X
\right\}\\
=&\left\{(v_1,\dots,v_{\ell})\in\mathbb{K}^{\ell}\,\middle|\,
\alpha_H\left(\sum_{i=1}^{\ell}v_ie_i\right)
=0\ {\rm for}\ H\in\mathscr{A}_X
\right\}\simeq X.\qquad\Box
\end{align*}

The linear map $\delta_X:V^{\ast}\rightarrow \mathbb{K}$ is defined by $\alpha\mapsto\delta_X(\alpha)$ for $\alpha\in V^{\ast}$.
Since $v_X\neq 0$, the map $\delta_X$ is not zero.
Then the image ${\rm Im}(\delta_X:V^{\ast}\rightarrow \mathbb{K})$ is $\mathbb{K}$.
We have the split exact sequence
\begin{align*}
\begin{CD}
0@>>>{\rm Ker}(\delta_X:V^{\ast}\rightarrow \mathbb{K})
@>>> V^{\ast}@>\delta_X>> \mathbb{K}@>>> 0
\end{CD}
\end{align*}
and there exists a section $\iota_X:\mathbb{K}\rightarrow V^{\ast}$ such that
$\delta_X\circ\iota_X=id_\mathbb{K}$, where $id_\mathbb{K}$ is the identity map of $\mathbb{K}$.
We define
\begin{align}
y_X:=\iota_X(1).
\end{align}
Let
\begin{align}
V^{\ast}_X:={\rm Ker}(\delta_X:V^{\ast}\rightarrow \mathbb{K})
\end{align}
and let $\{y_1,\dots,y_{\ell-1}\}$ be a $\mathbb{K}$-basis for $V^{\ast}_X$.
We note that $\delta_X(y_X)=1$
and $\delta_X(y_i)=0$ for any $1\leq i\leq \ell-1$.
Let $\{\delta_1,\dots,\delta_{\ell-1},\delta_{\ell}\}$ be the dual basis for
the basis $\{y_1,\dots,y_{\ell-1},y_X\}$ for
$V^{\ast}=V^{\ast}_X\oplus\iota_X(\mathbb{K})$. Then
\begin{align*}
\delta_X=\sum_{i=1}^{\ell-1}\delta_X(y_i)\delta_i
+\delta_X(y_X)\delta_{\ell}
=\delta_X(y_X)\delta_{\ell}=\delta_{\ell}.
\end{align*}
Hence we may consider
$\{y_1,\dots,y_{\ell-1},y_X\}$ as variables of $S$
with the partial derivatives $\{\delta_1,\dots,\delta_{\ell-1},\delta_X\}$.
\begin{figure}[t]
\begin{center}
\includegraphics[width=60mm]{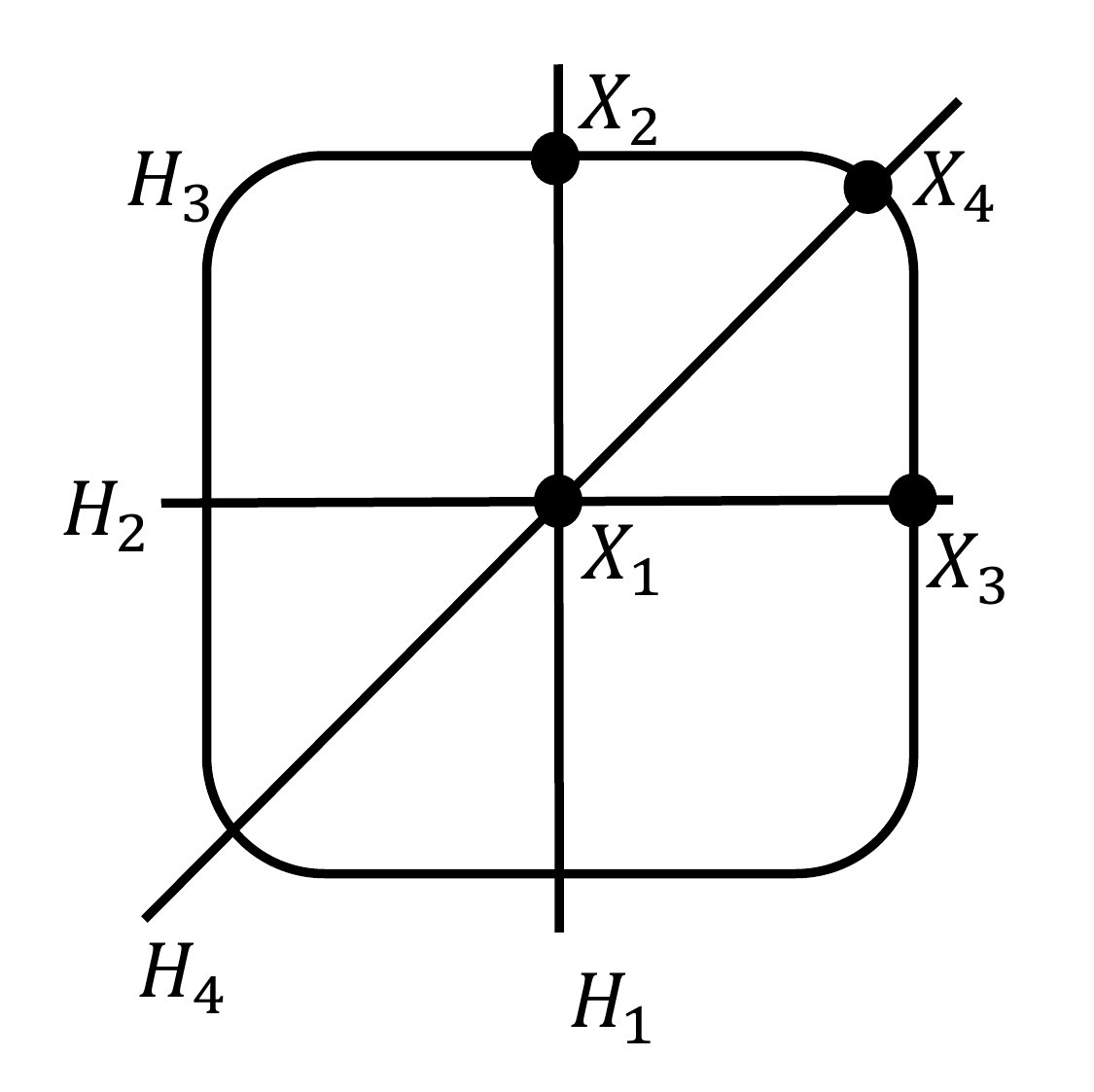}
\caption{The $3$-arrangement defined by $Q=x_1 x_2 x_3(x_1-x_2)$}
\label{fig-x1x2x3(x1-x2)}
\end{center}
\end{figure}
\begin{Examp}\label{ex-xyz(x-y)-delta}
Let $\ell=3$. Let us consider the $3$-arrangement $\mathscr{A}$
consisting of four hyperplanes
$H_1:=\{x_1=0\},H_2:=\{x_2=0\},H_3:=\{x_3=0\},H_4:=\{x_1-x_2=0\}$
(see Figure \ref{fig-x1x2x3(x1-x2)}).
The elements of $L(\mathscr{A})_1$ are
\begin{align*}
&X_1:=\{x_1=0,x_2=0\},\ X_2:=\{x_1=0,x_3=0\},\\
&X_3:=\{x_2=0,x_3=0\},\ X_4:=\{x_1=x_2,x_3=0\},
\end{align*}
and
\begin{align*}
\mathscr{A}_{X_1}=\{H_1,H_2,H_4\},\ \mathscr{A}_{X_2}=\{H_1,H_3\},\ 
\mathscr{A}_{X_3}=\{H_2,H_3\},\ \mathscr{A}_{X_4}=\{H_3,H_4\}.
\end{align*}
By Proposition \ref{prop-deltaX}, we may consider $\delta_X$ and $y_X$ $(X\in L(\mathscr{A})_1)$ as
\begin{align*}
&\delta_{X_1}=\partial_3,\ \delta_{X_2}=\partial_2,\ \delta_{X_3}=\partial_1,\ 
\delta_{X_4}=\partial_1+\partial_2,\\
&y_{X_1}=x_3,\ y_{X_2}=x_2,\ y_{X_3}=x_1,\ 
y_{X_4}=\frac{1}{2}(x_1+x_2).
\end{align*}
\end{Examp}

Let $S_X:={\rm Sym}(V^{\ast}_X)=\mathbb{K}[y_1,\dots,y_{\ell-1}]$
be the symmetric algebra of $V^{\ast}_X$.
Then $S=S_X\otimes_{\mathbb{K}}\mathbb{K}[y_X]$.
We note that $\alpha_H\in S_X$ if $H\in\mathscr{A}_X$
and $\alpha_H\not\in S_X$ if $H\in\mathscr{A}\setminus\mathscr{A}_X$.
%In addition we have $S_X\subsetneq S$ and
%$D^{(m)}(S_X)=\sum_{\bm{a}\in\mathbb{N}^{\ell-1},|\bm{a}|=m}
%S_X\delta^{\bm{a}}$.
We write
\begin{align}
Q_X:=Q(\mathscr{A}_X).
\end{align}
Let $V_X:=V/X$ be the quotient vector space of $V$ by $X$. Then $V^{\ast}_X$ is considered as the dual space of $V_X$.
Moreover, for an integer $j\geq 0$, the module
\begin{align*}
D^{(j)}(\mathscr{A}_X,V_X)=
\left\{\theta\in D^{(j)}(S_X)\,\middle|\,
\theta(Q_X S_X)\subseteq Q_XS_X\right\}
\end{align*}
is naturally embedded in $D^{(j)}(S)$ as an $S_X$-submodule
by extending its operators to $S=S_X\otimes_{\mathbb{K}}\mathbb{K}[y_X]$
so that their partial derivatives annihilate $y_X$.
We identify $D^{(j)}(\mathscr{A}_X,V_X)$ with its image under
this embedding and write $SD^{(j)}(\mathscr{A}_X,V_X)$
for the $S$-submodule of $D^{(j)}(S)$ generated by this image.

\section{Freeness when $\ell=3$ and $m\geq|\mathscr{A}|-2$}\label{sec-when-m>n-2}
In this section we prove that every $3$-arrangement $\mathscr{A}$ is $m$-free
when $m\geq n-2$.
An arrangement $\mathscr{A}$ is said to be essential
if $\bigcap_{H\in\mathscr{A}}H=\{0\}$.
We first investigate a nonessential arrangement, which is
a direct product of a $2$-arrangement and the empty $1$-arrangement.
Holm proved the following.
\begin{Prop}[Proposition I\hspace{-0.5mm}I\hspace{-0.5mm}I. 5.2 in \cite{Holm-phd}]\label{prop-2arrangements}
Let $(\mathscr{A},V)$ be a $2$-arrangement.
Then for any $m\geq 0$, $(\mathscr{A},V)$ is $m$-free with
\begin{align*}
{\rm exp}_m(\mathscr{A},V)=\left\{
\begin{array}{ll}
\{m,(|\mathscr{A}|-1)^{m}\}&(0\leq m\leq |\mathscr{A}|-1),\\
\{(|\mathscr{A}|-1)^{|\mathscr{A}|},|\mathscr{A}|^{m-|\mathscr{A}|+1}\}
&(m\geq |\mathscr{A}|),
\end{array}
\right.
\end{align*}
where $e^i\in{\rm exp}_m(\mathscr{A},V)$
means that the integer $e$ occurs $i$ times.
\end{Prop}
In addition Abe and the author characterised
the $m$-freeness of product arrangements as follows.
\begin{Them}[Theorem 2.1 in \cite{Abe-Nakashima}]\label{thm-AN-prod-arr}
Let $(\mathscr{A}_1,V_1)$ and $(\mathscr{A}_2,V_2)$ be arrangements
with $\dim(V_1)>0$ and $\dim(V_2)>0$.
The product arrangement $(\mathscr{A}_1\times\mathscr{A}_2,V_1\oplus V_2)$
is $m$-free if and only if both $(\mathscr{A}_1,V_1)$ and $(\mathscr{A}_2,V_2)$
are $i$-free for all $i\in\mathbb{N}$ with $1\leq i\leq m$.
Moreover, if $\exp_i(\mathscr{A}_1)=\{d_j^{(i)}\mid 1\leq j\leq s_i(\ell_1)\}$
and $\exp_{m-i}(\mathscr{A}_2)=\{e_k^{(m-i)}\mid 1\leq k\leq s_{m-i}(\ell_2)\}$
for $0\leq i\leq m$, then
$\exp_m(\mathscr{A}_1\times\mathscr{A}_2)=\bigcup_{i=0}^{m}
\{d_j^{(i)}+e_k^{(m-i)}\mid
1\leq j\leq s_i(\ell_1),\ 1\leq k\leq s_{m-i}(\ell_2)\}$.
\end{Them}
Proposition \ref{prop-2arrangements} and Theorem \ref{thm-AN-prod-arr}
imply that a nonessential $3$-arrangement is $m$-free for all $m\geq 0$.
If a nonessential $3$-arrangement $(\mathscr{A},V)$
is a product of a $2$-arrangement $(\mathscr{A}^{\prime},V^{\prime})$
and the empty $1$-arrangement, then ${\rm exp}_m(\mathscr{A},V)$ equals
$\bigcup_{j=0}^{m}{\rm exp}_j(\mathscr{A}^{\prime},V^{\prime})$.

We next investigate essential $3$-arrangements.
In the rest of this section,
we assume that $\mathscr{A}$ is an essential arrangement.
Then $L(\mathscr{A})_1$ is not empty.
We also assume that $\ell=3$ and $m\geq n-2$.
If $m>n-2$, then we take distinct hyperplanes $H_{n+1},\dots,H_{m+2}$ in $V$ with $H_i\not\in\mathscr{A}$ for $n+1\leq i\leq m+2$.
We define an arrangement $\widetilde{\mathscr{A}}$ by
\begin{align}
\widetilde{\mathscr{A}}:=\left\{
\begin{array}{ll}
\mathscr{A}\quad&{\rm if}\ m=n-2,\\
\mathscr{A}\cup\{H_{n+1},\dots,H_{m+2}\}\quad&{\rm if}\ m>n-2.
\end{array}
\right.
\end{align}
Then $\widetilde{\mathscr{A}}$ is also an essential arrangement and
$L(\widetilde{\mathscr{A}})_1$ is not empty. We write $\widetilde{n}:=|\widetilde{\mathscr{A}}|$ and then $m=\widetilde{n}-2$.
Let $X\in L(\widetilde{\mathscr{A}})_1$.
We use the variables $y_1,y_2,y_X$ and the derivatives
$\delta_1,\delta_2,\delta_X$ defined in Section \ref{sec-change-vari}.
Let
\begin{align}
i_X:=|\widetilde{\mathscr{A}}_X|-2\quad
{\rm and}\quad
\widetilde{P}_X:=
\prod_{H\in\widetilde{\mathscr{A}}\setminus\widetilde{\mathscr{A}}_X}
\alpha_H=\frac{Q(\widetilde{\mathscr{A}})}{Q(\widetilde{\mathscr{A}}_X)}.
\end{align}
Since $\dim(X)=1$, there exist at least two hyperplanes in $\widetilde{\mathscr{A}}_X$, that is, $|\widetilde{\mathscr{A}}_X|\geq 2$.
Then $i_X$ is a nonnegative integer.
We note that $\deg(\widetilde{P}_X)=\widetilde{n}-|\widetilde{\mathscr{A}}_X|=\widetilde{n}-2-|\widetilde{\mathscr{A}}_X|+2=m-i_X$.
For any $0\leq j\leq i_X$, let $B(X,j)$ be a $\mathbb{K}$-basis for the vector space $(S_X)_j$ of homogeneous polynomials of degree $j$ in $S_X$.
Then we define a set $B$ by
\begin{align}\label{eq-B}
B:=\bigcup_{X\in L(\widetilde{\mathscr{A}})_1}\bigcup_{j=0}^{i_X}
\left\{\widetilde{P}_X\,u\,y_X^{i_X-j}\,\middle|\,u\in B(X,j)\right\}
\subseteq S_m.
\end{align}
For any $0\leq j\leq i_X$,
let $\{\partial_{u}\mid u\in B(X,j)\}$ be the dual basis for $B(X,j)$,
that is, for $u,v\in B(X,j)$,
$\partial_{u}(v)=\left\{
\begin{array}{ll}
1\quad&{\rm if}\ u=v,\\
0\quad&{\rm if}\ u\neq v.
\end{array}
\right.$
We recall that $\delta_X(y_X)=1$ and $\deg(\widetilde{P}_X y_X^{i_X-j})=m-j$.
Since $\delta_X^{m-j}(\widetilde{P}_X y_X^{i_X-j})=
(m-j)!\prod_{H\in\widetilde{\mathscr{A}}\setminus\widetilde{\mathscr{A}}_X}
\delta_X(\alpha_H)\neq 0$
by Proposition \ref{prop-deltaX}, we can define a set $B^{\ast}$ by
\begin{align}
B^{\ast}:=\bigcup_{X\in L(\widetilde{\mathscr{A}})_1}\bigcup_{j=0}^{i_X}
\left\{\frac{1}{\delta_X^{m-j}(\widetilde{P}_X y_X^{i_X-j})}
\partial_{u}\delta_X^{m-j}\,\middle|\,u\in B(X,j)\right\}\subseteq
\sum_{|\bm{a}|=m}\mathbb{K}\partial^{\bm{a}}.
\end{align}

\begin{figure}[t]
\begin{center}
\includegraphics[width=60mm]{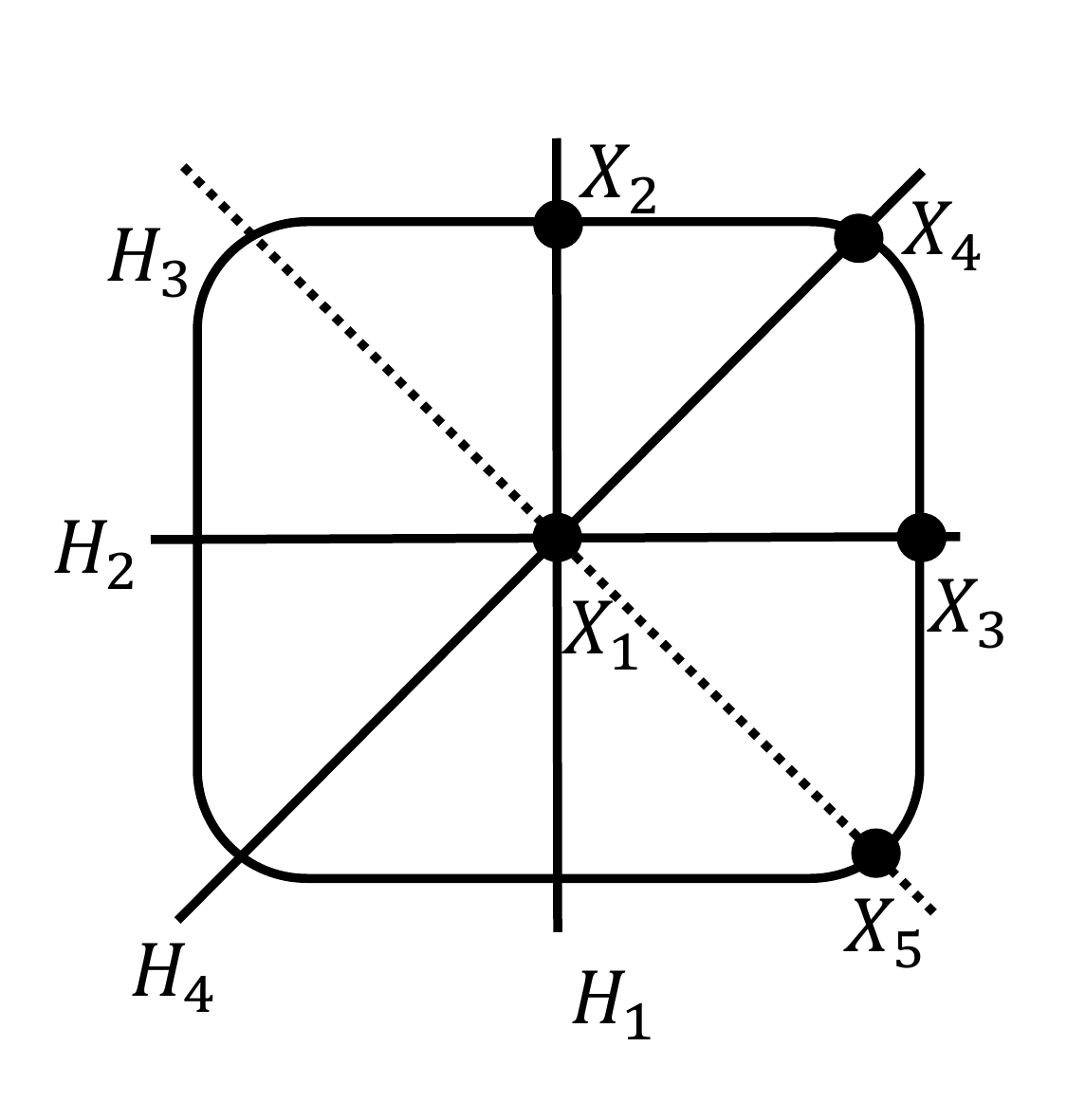}
\caption{$\widetilde{\mathscr{A}}$ in Example \ref{ex-xyz(x-y)-delta-PX} $(2)$}
\label{fig-Atilde}
\end{center}
\end{figure}

\begin{Examp}\label{ex-xyz(x-y)-delta-PX}
We consider the same $3$-arrangement $\mathscr{A}$
as Example \ref{ex-xyz(x-y)-delta}, i.e.,
$Q(\mathscr{A})=x_1 x_2 x_3(x_1-x_2)$.

$(1)$ Let $m=2=|\mathscr{A}|-2$.
In this case $\widetilde{\mathscr{A}}=\mathscr{A}$.
We recall that
$L(\mathscr{A})_1=\{X_1=\{x_1=0,x_2=0\},\ X_2=\{x_1=0,x_3=0\},
X_3=\{x_2=0,x_3=0\},\ X_4=\{x_1=x_2,x_3=0\}\}$.
Then $\delta_X$, $y_X$, $S_X$, $i_X$ and $\widetilde{P}_X$ $(X\in L(\mathscr{A})_1)$ are listed in the following table.
\begin{center}
\begin{tabular}{|c||c|c|c|c|}
\hline
 $L(\mathscr{A})_1$ &$X_1$&$X_2$&$X_3$&$X_4$\\
\hline
 $\delta_X$ & $\partial_3$ & $\partial_2$ & $\partial_1$ & $\partial_1+\partial_2$ \\
\hline
 $y_X$ & $x_3$ & $x_2$ & $x_1$ & $\frac{x_1+x_2}{2}$ \\
\hline
 $S_X$ & $\mathbb{K}[x_1,x_2]$ & $\mathbb{K}[x_1,x_3]$ & $\mathbb{K}[x_2,x_3]$ & $\mathbb{K}\left[\frac{x_1-x_2}{2},x_3\right]$ \\
\hline
 $i_X$ &$1$&$0$&$0$&$0$\\
\hline
 $\widetilde{P}_X$ &$x_3$&$x_2(x_1-x_2)$&$x_1(x_1-x_2)$&$x_1x_2$\\
\hline
\end{tabular}
\end{center}
We can take
$B(X_1,0)=\{1\}$,\ $B(X_1,1)=\{x_1,x_2\}$,\ $B(X_2,0)=\{1\}$,\ 
$B(X_3,0)=\{1\}$\ and\ $B(X_4,0)=\{1\}$.
Therefore
\begin{align*}
B&=\{x_3^2, x_3x_1, x_3x_2, x_2(x_1-x_2), x_1(x_1-x_2), x_1x_2\},\\
B^{\ast}&=\left\{\frac{1}{2}\partial_3^2, \partial_1\partial_3, 
\partial_2\partial_3, -\frac{1}{2}\partial_2^2, 
\frac{1}{2}\partial_1^2, \frac{1}{2}(\partial_1+\partial_2)^2\right\}.
\end{align*}

$(2)$
Let $m=3$. We take a hyperplane $H_5:=\{x_1+x_2=0\}$
and set $\widetilde{\mathscr{A}}=\mathscr{A}\cup\{H_5\}$.
The elements of $L(\widetilde{\mathscr{A}})_1$ are
$X_1:=\{x_1=0,x_2=0\}$, $X_2:=\{x_1=0,x_3=0\}$, $X_3:=\{x_2=0,x_3=0\}$,
$X_4:=\{x_1=x_2,x_3=0\}$, and $X_5:=\{x_1=-x_2,x_3=0\}$.
Then $\delta_X$, $y_X$, $S_X$, $i_X$, and $\widetilde{P}_X$ $(X\in L(\widetilde{\mathscr{A}})_1)$ are listed in the following table.
\begin{center}
\begin{tabular}{|c||c|c|c|c|c|}
\hline
$L(\widetilde{\mathscr{A}})_1$&$X_1$&$X_2$&$X_3$&$X_4$&$X_5$\\
\hline
 $\delta_X$ & $\partial_3$ & $\partial_2$ & $\partial_1$ & $\partial_1+\partial_2$ & $\partial_1-\partial_2$ \\
\hline
 $y_X$ & $x_3$ & $x_2$ & $x_1$ & $\frac{x_1+x_2}{2}$&$\frac{x_1-x_2}{2}$\\
\hline
 $i_X$ &$2$&$0$&$0$&$0$&$0$\\
\hline
 $S_X$ &$\mathbb{K}[x_1,x_2]$&$\mathbb{K}[x_1,x_3]$&$\mathbb{K}[x_2,x_3]$&$\mathbb{K}\left[\frac{x_1-x_2}{2},x_3\right]$&$\mathbb{K}\left[\frac{x_1+x_2}{2},x_3\right]$\\
\hline
 $\widetilde{P}_X$ &$x_3$&$x_2(x_1^2-x_2^2)$&$x_1(x_1^2-x_2^2)$&$x_1x_2(x_1+x_2)$&$x_1x_2(x_1-x_2)$\\
\hline
\end{tabular}
\end{center}
We can take
$B(X_1,0)=\{1\}$,\ $B(X_1,1)=\{x_1,x_2\}$,\ 
$B(X_1,2)=\{x_1^2,x_1x_2,x_2^2\}$,\ $B(X_2,0)=\{1\}$,\ 
$B(X_3,0)=\{1\}$,\ $B(X_4,0)=\{1\}$\, and\ $B(X_5,0)=\{1\}$.
\begin{align*}
B&=\left\{
\begin{array}{l}
x_3^3, x_3^2x_1, x_3^2x_2, x_3x_1^2, x_3x_1x_2, x_3x_2^2, x_2(x_1^2-x_2^2),\\ x_1(x_1^2-x_2^2), x_1x_2(x_1+x_2), x_1x_2(x_1-x_2)
\end{array}
\right\},\\
B^{\ast}&=\left\{
\begin{array}{l}
\frac{1}{6}\partial_3^3, \frac{1}{2}\partial_1\partial_3^2, \frac{1}{2}\partial_2\partial_3^2, \frac{1}{2}\partial_1^2\partial_3, \partial_1\partial_2\partial_3, \frac{1}{2}\partial_2^2\partial_3,\\
-\frac{1}{6}\partial_2^3, \frac{1}{6}\partial_1^3, \frac{1}{6}(\partial_1+\partial_2)^3, \frac{1}{6}(\partial_1-\partial_2)^3,
\end{array}
\right\}.
\end{align*}
\end{Examp}

We prepare Lemma \ref{lem-|A_XcapA_Y|<2} and Proposition \ref{prop-s_m=sum_Xs_iX} to prove that $B$ is a $\mathbb{K}$-basis for $S_m$ and that $B^{\ast}$ is a $\mathbb{K}$-basis for $\sum_{|\bm{a}|=m}\mathbb{K}\partial^{\bm{a}}$.
\begin{Lem}\label{lem-|A_XcapA_Y|<2}
Let $\ell=3$ and let $X,Y\in L(\widetilde{\mathscr{A}})_1$.
If $X\neq Y$, then $\delta_X^{m-i_X}(\widetilde{P}_Y f)=0$ for any $f\in S$
with $\deg(f)=i_Y$.
\end{Lem}

\noindent
{\it Proof.}
First note that $X\neq Y$ implies
$|\widetilde{\mathscr{A}}_X\cap\widetilde{\mathscr{A}}_Y|\leq 1$.
By Proposition \ref{prop-deltaX}, for any $f\in S$ with $\deg(f)=i_Y$,
\begin{align*}
\delta_X^{m-i_X}(\widetilde{P}_Y f)&=
\delta_X^{m-i_X}\left(\prod_{H\in\widetilde{\mathscr{A}}
\setminus\widetilde{\mathscr{A}}_Y}\alpha_Hf\right)\\
&=\prod_{H\in\widetilde{\mathscr{A}}_X\setminus
(\widetilde{\mathscr{A}}_X\cap\widetilde{\mathscr{A}}_Y)}\alpha_H
\cdot\delta_X^{m-i_X}\left(\prod_{H\in\widetilde{\mathscr{A}}
\setminus(\widetilde{\mathscr{A}}_X\cup\widetilde{\mathscr{A}}_Y)}
\alpha_H f\right).
\end{align*}
Here if $X\neq Y$, then
\begin{align*}
\deg\left(\prod_{H\in\widetilde{\mathscr{A}}
\setminus(\widetilde{\mathscr{A}}_X\cup\widetilde{\mathscr{A}}_Y)}
\alpha_H f\right)
&=|\widetilde{\mathscr{A}}\ |-
|\widetilde{\mathscr{A}}_X\cup\widetilde{\mathscr{A}}_Y|+i_Y\\
&=|\widetilde{\mathscr{A}}\ |-\left(|\widetilde{\mathscr{A}}_X|+
|\widetilde{\mathscr{A}}_Y|
-|\widetilde{\mathscr{A}}_X\cap\widetilde{\mathscr{A}}_Y|\right)
+\left(|\widetilde{\mathscr{A}}_Y|-2\right)\\
&=|\widetilde{\mathscr{A}}\ |-|\widetilde{\mathscr{A}}_X|
+|\widetilde{\mathscr{A}}_X\cap\widetilde{\mathscr{A}}_Y|-2\\
&<|\widetilde{\mathscr{A}}\ |-|\widetilde{\mathscr{A}}_X|=m-i_X.
\end{align*}
Hence if $X\neq Y$, then $\delta_X^{m-i_X}(\widetilde{P}_Y f)=0$.\hfill$\Box$

\begin{Prop}\label{prop-s_m=sum_Xs_iX}
Let $\mathscr{A}$ be an essential $3$-arrangement.
\begin{enumerate}
\item[$(1)$] We have the set equation
\begin{align}\label{eq-num-of-subsets-with|B|=l-1}
\left\{\mathscr{B}\subseteq \widetilde{\mathscr{A}}\ \middle|\,
|\mathscr{B}|=2\right\}
=\bigsqcup_{X\in L(\widetilde{\mathscr{A}})_1}
\left\{\mathscr{B}\subseteq \widetilde{\mathscr{A}}_X\,\middle|\,
|\mathscr{B}|=2\right\},
\end{align}
where the right hand side is a disjoint union.
\item[$(2)$]
\begin{align*}
s_m(3)=\sum_{X\in L(\widetilde{\mathscr{A}})_1}s_{i_X}(3).
\end{align*}
\end{enumerate}
\end{Prop}

\noindent
{\it Proof.}
$(1)$\ It is obvious that
$\left\{\mathscr{B}\subseteq \widetilde{\mathscr{A}}_X\,\middle|\,
|\mathscr{B}|=2\right\}\subseteq
\left\{\mathscr{B}\subseteq \widetilde{\mathscr{A}}\ \middle|\,
|\mathscr{B}|=2\right\}$
for any $X\in L(\widetilde{\mathscr{A}})_1$.
Conversely let
$\mathscr{B}=\{H,H^{\prime}\}\subseteq \widetilde{\mathscr{A}}$.
We set $X:=H\cap H^{\prime}\in L(\widetilde{\mathscr{A}})_1$ and then $\mathscr{B}\subseteq\widetilde{\mathscr{A}}_X$.
Therefore $\mathscr{B}$ belongs to the right hand side of \eqref{eq-num-of-subsets-with|B|=l-1}.

Next we prove that if $X,Y\in L(\widetilde{\mathscr{A}})_1$ with $X\neq Y$, then $\Omega:=\{\mathscr{B}\subseteq \widetilde{\mathscr{A}}_X\mid|\mathscr{B}|=2\}\cap\{\mathscr{B}\subseteq \widetilde{\mathscr{A}}_Y\mid|\mathscr{B}|=2\}$ is the empty set.
We suppose that $\Omega$ is not empty
and let $\mathscr{B}=\{H,H^{\prime}\}\in \Omega$.
Since $\dim(H\cap H^{\prime})=1$ and $X\subseteq H\cap H^{\prime}$,
we have $X=H\cap H^{\prime}$.
Similarly $Y=H\cap H^{\prime}$ and hence $X=Y$.
The argument above imply that the right hand side of
\eqref{eq-num-of-subsets-with|B|=l-1} is a disjoint union.

$(2)$\ We have
\begin{align*}
s_m(3)&=\binom{m+2}{2}
=\binom{\widetilde{n}}{2}
=\sum_{X\in L(\widetilde{\mathscr{A}})_1}
\binom{|\widetilde{\mathscr{A}}_X|}{2}
=\sum_{X\in L(\widetilde{\mathscr{A}})_1}\binom{i_X+2}{2}
=\sum_{X\in L(\widetilde{\mathscr{A}})_1}s_{i_X}(3)
\end{align*}
by the equation \eqref{eq-num-of-subsets-with|B|=l-1}.
\hfill$\Box$

\begin{Prop}\label{prop-B-basis-S_M}
Let $\ell=3$. The set $B$ is a $\mathbb{K}$-basis for $S_m$
and $B^{\ast}$ is a $\mathbb{K}$-basis for
$\sum_{|\bm{a}|=m}\mathbb{K}\partial^{\bm{a}}$.
\end{Prop}

\noindent
{\it Proof.}
For $X\in L(\widetilde{\mathscr{A}})_1$,
$0\leq j\leq i_X$, and $u\in B(X,j)$, let
\begin{align*}
b_{X,j,u}
&:=\widetilde{P}_X u y_X^{i_X-j},\\
c_{X,j}
&:=\delta_X^{m-j}
  (\widetilde{P}_X y_X^{i_X-j})\in\mathbb{K}\setminus\{0\},\\
\lambda_{X,j,u}
&:=\frac{1}{c_{X,j}}\partial_u\delta_X^{m-j}.
\end{align*}
We consider the matrix
\begin{align*}
M:=
\left(
\lambda_{X,j,u}(b_{Y,k,v})
\right)_{(X,j,u),(Y,k,v)},
\end{align*}
whose rows and columns are indexed by
\begin{align*}
\left\{
(X,j,u)\,\middle|\,
X\in L(\widetilde{\mathscr{A}})_1,\ 
0\leq j\leq i_X,\ 
u\in B(X,j)
\right\}.
\end{align*}

First suppose that $X\neq Y$.
Since $\deg(vy_Y^{i_Y-k})=i_Y$, Lemma
\ref{lem-|A_XcapA_Y|<2} implies that
\begin{align*}
\delta_X^{m-j}(b_{Y,k,v})
&=
\delta_X^{i_X-j}\delta_X^{m-i_X}
\left(\widetilde{P}_Y v y_Y^{i_Y-k}\right)
=0.
\end{align*}
Hence $\lambda_{X,j,u}(b_{Y,k,v})=0$.

Next suppose that $X=Y$.
Since $v\in S_X$, we have $\delta_X(v)=0$ and thus
\begin{align*}
\delta_X^{m-j}(b_{X,k,v})
=
v\,\delta_X^{m-j}
\left(\widetilde{P}_X y_X^{i_X-k}\right).
\end{align*}
If $j<k$, then
$\deg(\widetilde{P}_X y_X^{i_X-k})=m-k<m-j$.
Therefore
$\delta_X^{m-j}
(\widetilde{P}_X y_X^{i_X-k})=0$, and hence
$\lambda_{X,j,u}(b_{X,k,v})=0$.

If $j=k$, then
$\delta_X^{m-j}
(\widetilde{P}_X y_X^{i_X-j})=c_{X,j}$
is a nonzero constant. It follows that
\begin{align}\label{eq-bualbasis}
\lambda_{X,j,u}(b_{X,j,v})
=\partial_u(v)
=\left\{
\begin{array}{ll}
1 & {\rm if}\ u=v,\\
0 & {\rm if}\ u\neq v.
\end{array}
\right.
\end{align}

Order the rows and columns of $M$ first by $X$ and,
for each fixed $X$, by increasing $j$,
using the same ordering of $B(X,j)$ for the rows and columns.
The preceding calculations show that $M$ is block diagonal
with respect to $X$, and that each diagonal block is
block lower triangular with identity diagonal blocks.
Therefore $\det(M)=1$.

Consequently, the indexed families
$\{b_{X,j,u}\}$ and $\{\lambda_{X,j,u}\}$
are both linearly independent over $\mathbb{K}$.
In particular, neither family has repeated elements,
so these families give the sets $B$ and $B^{\ast}$,
respectively.
By Proposition \ref{prop-s_m=sum_Xs_iX} $(2)$,
their common cardinality is
\begin{align*}
\sum_{X\in L(\widetilde{\mathscr{A}})_1}
\sum_{j=0}^{i_X}|B(X,j)|
&=
\sum_{X\in L(\widetilde{\mathscr{A}})_1}
\sum_{j=0}^{i_X}s_j(2)\\
&=
\sum_{X\in L(\widetilde{\mathscr{A}})_1}s_{i_X}(3)
=s_m(3).
\end{align*}
Since
\begin{align*}
\dim_{\mathbb{K}}S_m
=
\dim_{\mathbb{K}}
\left(\sum_{|\bm{a}|=m}\mathbb{K}\partial^{\bm{a}}\right)
=s_m(3),
\end{align*}
$B$ and $B^{\ast}$ are bases for the respective vector spaces.
\hfill$\Box$
\vspace{2mm}

Let
\begin{align}
P_X:=\prod_{H\in\mathscr{A}\setminus\mathscr{A}_X}\alpha_H=\frac{Q}{Q_X}.
\end{align}
Then we have
\begin{align*}
\widetilde{P}_X=
\frac{Q(\widetilde{\mathscr{A}}\ )}{Q(\widetilde{\mathscr{A}}_X)}
=\frac{Q\prod_{H\in\widetilde{\mathscr{A}}\setminus\mathscr{A}}\alpha_H}{Q_X\prod_{H\in\widetilde{\mathscr{A}}_X\setminus\mathscr{A}_X}\alpha_H}
=P_X\prod_{H\in\widetilde{\mathscr{A}}\setminus(\mathscr{A}\cup\widetilde{\mathscr{A}}_X)}\alpha_H.
\end{align*}
Holm \cite{Holm-phd,Holm} proved the following (see also \cite[Corollary 2.5]{Abe-Nakashima}).
\begin{Prop}[Holm]\label{prop-check-opes}
Let $(\mathscr{A},V)$ be an $\ell$-arrangement. Then
\begin{align*}
D^{(m)}(\mathscr{A},V)=\bigcap_{H\in\mathscr{A}}
\left\{\theta\in D^{(m)}(S)\,\middle|\,
\theta(\alpha_{H}f)\in \alpha_{H}S\ 
{\rm for\ all}\ f\in S_{m-1}
\right\}.
\end{align*}
\end{Prop}

For operators $\theta_{1},\dots,\theta_{s_m(\ell)}\in D^{(m)}(S)$,
a coefficient matrix $M_{m}(\theta_{1},\dots,\theta_{s_m(\ell)})$
is an $s_m(\ell)\times s_m(\ell)$ matrix defined by
$$M_{m}(\theta_{1},\dots,\theta_{s_m(\ell)}):=
\left(\frac{\theta_i\left(x^{\bm{a}}\right)}{\bm{a}!}
\right)_{1\leq i\leq s_m(\ell),\ \bm{a}\in\{\bm{b}\in\mathbb{N}^{\ell}\mid|\bm{b}|=m\}}.$$
The following criterion for $D^{(m)}(\mathscr{A},V)$ is first given by Saito \cite{Saito-criterion} for $D^{(1)}(\mathscr{A},V)$ and is generalized by Holm \cite{Holm-phd} for $D^{(m)}(\mathscr{A},V)$ (see also Theorem 3.8 and Theorem 3.9 in \cite{Abe-Nakashima}).
\begin{Them}[Saito's criterion]\label{thm-saito's-criterion}
Let $\mathscr{A}$ be an $\ell$-arrangement and
let $\theta_{1},\dots,\theta_{s_m(\ell)}\in D^{(m)}(\mathscr{A},V)$ be homogeneous operators.
Then we have the following:
\begin{itemize}
\item[$(1)$]\ $\det M_{m}(\theta_{1},\dots,\theta_{s_m(\ell)})
=cQ^{s_{m-1}(\ell)}$ for some $c\in \mathbb{K}\setminus\{0\}$
if and only if the set $\{\theta_{1},\dots,\theta_{s_m(\ell)}\}$
is a free basis for $D^{(m)}(\mathscr{A})$ over $S$.
\item[$(2)$]\ $\mathscr{A}$ is $m$-free with homogeneous basis
$\theta_{1},\dots,\theta_{s_m(\ell)}$ if and only if these operators are linearly independent over $S$ and $\sum_{j=1}^{s_{m}(\ell)}\deg\left(\theta_{j}\right)=s_{m-1}(\ell)|\mathscr{A}|$.
\end{itemize}
\end{Them}

We use Proposition \ref{prop-B-basis-S_M} and Saito's criterion
to prove Theorem \ref{thm-D^m(A)=oplus_XD^i_Xdelta_X^m-i_X},
which is the key to the proof of our main theorem.

\begin{Them}\label{thm-D^m(A)=oplus_XD^i_Xdelta_X^m-i_X}
Let $\mathscr{A}$ be an essential $3$-arrangement.
If $m\geq n-2$, then $\mathscr{A}$ is $m$-free and
\begin{align}\label{eq-D^m(A)=oplus_XD^i_Xdelta_X^m-i_X}
D^{(m)}(\mathscr{A},V)
=
\bigoplus_{X\in L(\widetilde{\mathscr{A}})_1}
\bigoplus_{j=0}^{i_X}
SP_XD^{(j)}(\mathscr{A}_X,V_X)\delta_X^{m-j}
\end{align}
as $S$-modules.
\end{Them}

\noindent
{\it Proof.}
Let $X\in L(\widetilde{\mathscr{A}})_1$,
let $0\leq j\leq i_X$, and let
$\theta\in D^{(j)}(\mathscr{A}_X,V_X)$.
Since $\delta_X(Q_X)=0$ and $\theta$ preserves $Q_XS$,
for any $f\in S$ we have
\begin{align*}
P_X\theta\delta_X^{m-j}(Qf)
&=
P_X\theta\delta_X^{m-j}(Q_XP_Xf)\\
&=
P_X\theta\left(Q_X\delta_X^{m-j}(P_Xf)\right)
\in P_XQ_XS=QS.
\end{align*}
Therefore
$P_X\theta\delta_X^{m-j}\in D^{(m)}(\mathscr{A},V)$.

By Proposition \ref{prop-2arrangements}, for each $X$ and $j$
we can take a homogeneous free $S_X$-basis
$\eta_{X,j,1},\dots,\eta_{X,j,j+1}$
for $D^{(j)}(\mathscr{A}_X,V_X)$ such that
\begin{align*}
\sum_{r=1}^{j+1}\deg(\eta_{X,j,r})=j|\mathscr{A}_X|.
\end{align*}
For $j=0$, we take $\eta_{X,0,1}=1$.
Let
\begin{align*}
\mathcal{G}:=
\left\{
P_X\eta_{X,j,r}\delta_X^{m-j}
\,\middle|\,
X\in L(\widetilde{\mathscr{A}})_1,\ 
0\leq j\leq i_X,\ 
1\leq r\leq j+1
\right\}.
\end{align*}
All these operators are homogeneous and belong to
$D^{(m)}(\mathscr{A},V)$.

By Proposition \ref{prop-B-basis-S_M}, the set
\begin{align*}
\mathcal{E}:=
\bigcup_{X\in L(\widetilde{\mathscr{A}})_1}
\bigcup_{j=0}^{i_X}
\left\{
\partial_u\delta_X^{m-j}
\,\middle|\,u\in B(X,j)
\right\}
\end{align*}
is a $\mathbb{K}$-basis for
$\sum_{|\bm{a}|=m}\mathbb{K}\partial^{\bm{a}}$,
since its elements are nonzero constant multiples
of the elements of $B^{\ast}$.
Consequently,
\begin{align*}
D^{(m)}(S)
=
\bigoplus_{X\in L(\widetilde{\mathscr{A}})_1}
\bigoplus_{j=0}^{i_X}
\bigoplus_{u\in B(X,j)}
S\partial_u\delta_X^{m-j}.
\end{align*}
For each fixed $X$ and $j$, the corresponding elements
of $\mathcal{G}$ lie in the summand indexed by $(X,j)$.
They are linearly independent over $S$, since
$\eta_{X,j,1},\dots,\eta_{X,j,j+1}$ are linearly independent
over $S_X$, $S=S_X[y_X]$, and $P_X\neq 0$.
Thus $\mathcal{G}$ is linearly independent over $S$.
Moreover, Proposition \ref{prop-s_m=sum_Xs_iX} $(2)$ gives
\begin{align*}
|\mathcal{G}|
=
\sum_{X\in L(\widetilde{\mathscr{A}})_1}
\sum_{j=0}^{i_X}(j+1)
=
\sum_{X\in L(\widetilde{\mathscr{A}})_1}s_{i_X}(3)
=s_m(3).
\end{align*}

Since $\deg(P_X)=n-|\mathscr{A}_X|$, the sum of the degrees of
the elements of $\mathcal{G}$ is
\begin{align*}
\sum_{g\in\mathcal{G}}\deg(g)
&=
\sum_{X\in L(\widetilde{\mathscr{A}})_1}
\sum_{j=0}^{i_X}
\left((j+1)(n-|\mathscr{A}_X|)+j|\mathscr{A}_X|\right)\\
&=
ns_m(3)-
\sum_{X\in L(\widetilde{\mathscr{A}})_1}|\mathscr{A}_X|(i_X+1).
\end{align*}
For each $H\in\mathscr{A}$, the same argument as in the proof of
Proposition \ref{prop-s_m=sum_Xs_iX} $(1)$ shows that
each $H'\in\widetilde{\mathscr{A}}\setminus\{H\}$
belongs to $\widetilde{\mathscr{A}}_X\setminus\{H\}$
for exactly one $X\in L(\widetilde{\mathscr{A}})_1$
with $X\subseteq H$, namely $X=H\cap H'$.
Thus
\begin{align*}
\widetilde{\mathscr{A}}\setminus\{H\}
=
\bigsqcup_{\substack{X\in L(\widetilde{\mathscr{A}})_1\\
X\subseteq H}}
\left(\widetilde{\mathscr{A}}_X\setminus\{H\}\right).
\end{align*}
Hence
\begin{align*}
\sum_{\substack{X\in L(\widetilde{\mathscr{A}})_1\\
X\subseteq H}}(i_X+1)
=
\sum_{\substack{X\in L(\widetilde{\mathscr{A}})_1\\
X\subseteq H}}
\left(|\widetilde{\mathscr{A}}_X|-1\right)
=
|\widetilde{\mathscr{A}}|-1=m+1.
\end{align*}
Therefore
\begin{align*}
\sum_{X\in L(\widetilde{\mathscr{A}})_1}|\mathscr{A}_X|(i_X+1)
=\sum_{X\in L(\widetilde{\mathscr{A}})_1}\sum_{H\in\mathscr{A}_X}(i_X+1)
=\sum_{H\in\mathscr{A}}
\sum_{\substack{X\in L(\widetilde{\mathscr{A}})_1\\
X\subseteq H}}(i_X+1)
=n(m+1).
\end{align*}
It follows that
\begin{align*}
\sum_{g\in\mathcal{G}}\deg(g)
=
n\left(s_m(3)-(m+1)\right)
=ns_{m-1}(3).
\end{align*}
By Theorem \ref{thm-saito's-criterion} $(2)$,
$\mathcal{G}$ is a free $S$-basis for
$D^{(m)}(\mathscr{A},V)$.
For each $X$ and $j$, the corresponding elements of
$\mathcal{G}$ form an $S$-basis for
$SP_XD^{(j)}(\mathscr{A}_X,V_X)\delta_X^{m-j}$.
Thus these submodules form the direct sum in
\eqref{eq-D^m(A)=oplus_XD^i_Xdelta_X^m-i_X}.
\hfill$\Box$
\vspace{2mm}

\noindent
{\it Remark.}
Let $\mathscr{A}$ be a generic $\ell$-arrangement.
It is known that $\mathscr{A}$ is $m$-free if and only if $m\geq n-\ell+1$ \cite{Holm-phd,NOS}.
The ``if'' part of this result coincides with Theorem \ref{thm-D^m(A)=oplus_XD^i_Xdelta_X^m-i_X} when $\ell=3$.
In addition if $\ell=3$ and $n=4$, then $\mathscr{A}$ is not $1$-free.
This is an example that a $3$-arrangement is not $m$-free when $m=n-3$.
\vspace{2mm}

\begin{Examp}\label{ex-basis-D^mA}
We consider the same $3$-arrangement $\mathscr{A}$
as Example \ref{ex-xyz(x-y)-delta}
and Example \ref{ex-xyz(x-y)-delta-PX}, i.e.,
$Q(\mathscr{A})=x_1 x_2 x_3(x_1-x_2)$.

$(1)$ Let $m=2$. Then $\widetilde{\mathscr{A}}=\mathscr{A}$
and $\widetilde{P}_X=P_X$ for any $X\in L(\widetilde{\mathscr{A}})_1$.
By Saito's criterion, $\left\{x_1\partial_1+x_2\partial_2,x_2(x_1-x_2)\partial_2\right\}$ is a free basis for $D^{(1)}(\mathscr{A}_{X_1},V_{X_1})$.
Hence by Theorem \ref{thm-D^m(A)=oplus_XD^i_Xdelta_X^m-i_X},
\begin{align*}
&D^{(2)}(\mathscr{A},V)\\
=&SP_{X_1}D^{(0)}(\mathscr{A}_{X_1},V_{X_1})\delta_{X_1}^2
\oplus SP_{X_1}D^{(1)}(\mathscr{A}_{X_1},V_{X_1})\delta_{X_1}^1
\oplus SP_{X_2}D^{(0)}(\mathscr{A}_{X_2},V_{X_2})\delta_{X_2}^2\\
&\oplus SP_{X_3}D^{(0)}(\mathscr{A}_{X_3},V_{X_3})\delta_{X_3}^2
\oplus SP_{X_4}D^{(0)}(\mathscr{A}_{X_4},V_{X_4})\delta_{X_4}^2\\
=&Sx_3\partial_3^2\oplus Sx_3(x_1\partial_1+x_2\partial_2)\partial_3
\oplus Sx_3 x_2(x_1-x_2)\partial_2\partial_3\oplus Sx_2(x_1-x_2)\partial_2^2\\
&\oplus Sx_1(x_1-x_2)\partial_1^2\oplus Sx_1 x_2(\partial_1+\partial_2)^2
\end{align*}
and ${\rm exp}_{2}(\mathscr{A},V)=\{1,2,3,2,2,2\}=\{1^1,2^4,3^1\}$.

$(2)$
Let $m=3$, let $H_5=\{x_1+x_2=0\}$, and let $\widetilde{\mathscr{A}}=\mathscr{A}\cup\{H_5\}$.
Then $P_X$ $(X\in L(\widetilde{\mathscr{A}}))_1$
are given by
\begin{align*}
P_{X_1}=x_3,\ P_{X_2}=x_2(x_1-x_2),\ P_{X_3}=x_1(x_1-x_2),\ 
P_{X_4}=x_1x_2,\ P_{X_5}=x_1x_2(x_1-x_2).
\end{align*}
By Saito's criterion and Theorem \ref{thm-D^m(A)=oplus_XD^i_Xdelta_X^m-i_X},
we have
\begin{align*}
D^{(3)}(\mathscr{A},V)
=&Sx_3\partial_3^3\oplus Sx_3(x_1\partial_1+x_2\partial_2)\partial_3^2
\oplus Sx_3 x_2(x_1-x_2)\partial_2\partial_3^2
\oplus Sx_3x_1(x_1-x_2)\partial_1^2\partial_3\\
&\oplus Sx_3x_2(x_1-x_2)\partial_2^2\partial_3
\oplus Sx_3x_1 x_2(\partial_1+\partial_2)^2\partial_3
\oplus Sx_2(x_1-x_2)\partial_2^3\\
&\oplus Sx_1(x_1-x_2)\partial_1^3\oplus Sx_1 x_2(\partial_1+\partial_2)^3
\oplus Sx_1 x_2(x_1-x_2)(\partial_1-\partial_2)^3
\end{align*}
and ${\rm exp}_{3}(\mathscr{A},V)=\{1,2,3,3,3,3,2,2,2,3\}=\{1^1,2^4,3^5\}$.
\end{Examp}

We remark that the multi-set of $m$-exponents does not change if we take other hyperplanes $H^{\prime}_{n+1},\dots,H^{\prime}_{m+2}$ instead of $H_{n+1},\dots,H_{m+2}$.

\begin{figure}[t]
\begin{center}
\includegraphics[width=60mm]{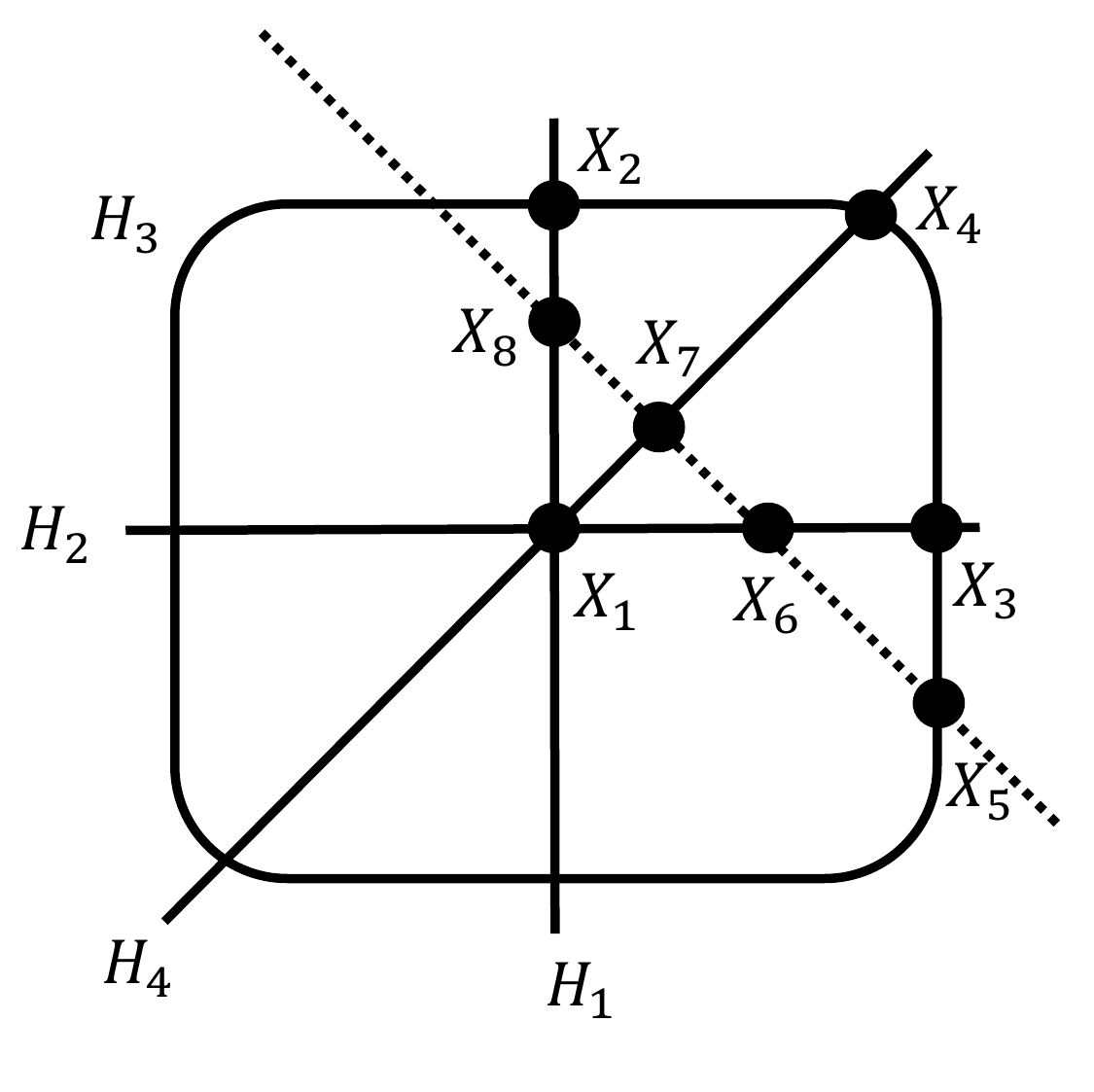}
\caption{$\widetilde{\mathscr{A}}$ in Example \ref{ex-anotherAtilde}}
\label{fig-anothorAtilde}
\end{center}
\end{figure}
\begin{Examp}\label{ex-anotherAtilde}
We consider the same $3$-arrangement $\mathscr{A}$
as Example \ref{ex-xyz(x-y)-delta}, i.e.,
$Q(\mathscr{A})=x_1 x_2 x_3(x_1-x_2)$.
Let $m=3$. In Example \ref{ex-basis-D^mA} we have already constructed
a free $S$-basis for $D^{(3)}(\mathscr{A},V)$.
In this example we construct another $S$-basis for $D^{(3)}(\mathscr{A},V)$
by taking a hyperplane $H^{\prime}_5:=\{x_1+x_2-x_3=0\}$.
We set $\widetilde{\mathscr{A}^{\prime}}=\mathscr{A}\cup\{H^{\prime}_5\}$.
Then the elements of $L(\widetilde{\mathscr{A}^{\prime}})_1$ are
$X_1:=\{x_1=0,x_2=0\}$, $X_2:=\{x_1=0,x_3=0\}$, $X_3:=\{x_2=0,x_3=0\}$,
$X_4:=\{x_1=x_2,x_3=0\}$, $X_5:=\{x_1=-x_2,x_3=0\}$, $X_6:=\{x_1=x_3,x_2=0\}$,
$X_7:=\{x_1=x_2=\frac{1}{2}x_3\}$, and $X_8:=\{x_2=x_3,x_1=0\}$.
Meanwhile $\delta_X$, $i_X$, and $P_X$ are listed in the following table.
\begin{flushleft}
\begin{tabular}{|c||c|c|c|c|c|}
\hline
$L(\widetilde{\mathscr{A}^{\prime}})_1$&$X_1$&$X_2$&$X_3$&$X_4$&$X_5$\\
\hline
 $\delta_X$ & $\partial_3$ & $\partial_2$ & $\partial_1$ & $\partial_1+\partial_2$ & $\partial_1-\partial_2$ \\
\hline
 $i_X$ &$1$&$0$&$0$&$0$&$0$\\
\hline
 $P_X$ & \qquad$x_3$\qquad\ &$x_2(x_1-x_2)$&$x_1(x_1-x_2)$&\ \qquad$x_1x_2$\qquad\ &$x_1x_2(x_1-x_2)$\quad\\
\hline
\end{tabular}\\
\vspace{2mm}
\begin{tabular}{|c||c|c|c|}
\hline
$L(\widetilde{\mathscr{A}^{\prime}})_1$&$X_6$&$X_7$&$X_8$\\
\hline
 $\delta_X$ & $\partial_1+\partial_3$ & $\partial_1+\partial_2+2\partial_3$ & $\partial_2+\partial_3$ \\
\hline
 $i_X$ &$0$&$0$&$0$\\
\hline
 $P_X$ &$x_1x_3(x_1-x_2)$&$x_1x_2x_3$&$x_2x_3(x_1-x_2)$\\
\hline
\end{tabular}
\end{flushleft}
By Theorem \ref{thm-D^m(A)=oplus_XD^i_Xdelta_X^m-i_X} and
Saito's criterion, we have
\begin{align*}
D^{(3)}(\mathscr{A},V)
=&SP_{X_1}\delta_{X_1}^3\oplus SP_{X_1}(x_1\partial_1+x_2\partial_2)
\delta_{X_1}^2\oplus SP_{X_1}x_2(x_1-x_2)\partial_2\delta_{X_1}^2
\oplus \bigoplus_{j=2}^{8}SP_{X_{j}}\delta_{X_{j}}^3
\end{align*}
and ${\rm exp}_{3}(\mathscr{A},V)=\{1,2,3,2,2,2,3,3,3,3\}=\{1^1,2^4,3^5\}$.
The multi-set of $m$-exponents coincides with that of Example \ref{ex-basis-D^mA}.
\end{Examp}

\section{Exponents}\label{sec-exp-3arr}
In this section we determine the $m$-exponents of a $3$-arrangement $\mathscr{A}$ for $m\geq n-2$, which depend only on the poset $L(\mathscr{A})$.
The following lemma holds for any $\ell\geq 3$.
\begin{Lem}\label{lem-addH_0}
Let $\mathscr{A}$ be an essential arrangement.
There exist a hyperplane $H_0\subseteq V$
with $H_0\not\in\mathscr{A}$ such that for any $X\in L(\mathscr{A})_1$,
$H_{0}\cap X=\{0\}$.
\end{Lem}

\noindent
{\it Proof.}
Let $H_0=\{a_1x_1+\cdots+a_{\ell}x_{\ell}=0\}$ be a hyperplane in $V$.
Let $X\in L(\mathscr{A})_1$ and let $0\neq v_X\in X$.
Then $X=\mathbb{K}v_X$. 
We note that $v_X=\left(x_1(v_X),\dots,x_{\ell}(v_X)\right)$. Here
\begin{align*}
H_0\cap X=\{0\}\ \Leftrightarrow\ v_X\not\in H_0
\ \Leftrightarrow\ a_1x_1(v_X)+\cdots+a_{\ell}x_{\ell}(v_X)\neq 0.
\end{align*}
Therefore $H_{0}\cap X=\{0\}$ for all $X\in L(\mathscr{A})_1$
$\Leftrightarrow$ $(a_1,\dots,a_{\ell})$ belongs to the complement of
$\bigcup_{X\in L(\mathscr{A})_1}
\{(x_1,\dots,x_{\ell})\mid x_1(v_X)x_1+\cdots+x_{\ell}(v_X)x_{\ell}=0\}$.
Since $\mathbb{K}$ is an infinite field,
we can take $(a_1,\dots,a_{\ell})$ from the complement of
$\bigcup_{X\in L(\mathscr{A})_1}
\{(x_1,\dots,x_{\ell})\mid x_1(v_X)x_1+\cdots+x_{\ell}(v_X)x_{\ell}=0\}$.
Then $H_0=\{a_1x_1+\cdots+a_{\ell}x_{\ell}=0\}$ is a desired hyperplane.
\hfill$\Box$
\vspace{2mm}

Let $\ell=3$, let $m\geq n-2$ and let $\mathscr{A}$ be essential.
The argument in Section \ref{sec-when-m>n-2} does not depend on a choice of hyperplanes $H_{n+1},\dots,H_{m+2}$.
By Lemma \ref{lem-addH_0}, if $m>n-2$ then there exist hyperplanes $H^{\prime}_{n+1},\dots,H^{\prime}_{m+2}$ in $V$ such that the next condition (A) is satisfied.
\begin{enumerate}
\item[(A)] For any $n+1\leq i\leq m+2$ and
for any $X\in L(\mathscr{A}\cup\{H^{\prime}_{n+1},\dots,H^{\prime}_{i-1}\})_1$,
$H^{\prime}_i\cap X=\{0\}$, where
$\{H^{\prime}_{n+1},\dots,H^{\prime}_{i-1}\}=\emptyset$
if $i=n+1$.
\end{enumerate}
We define $\mathscr{G}:=\mathscr{A}$ if $m=n-2$
and $\mathscr{G}:=\mathscr{A}\cup
\{H^{\prime}_{n+1},\dots,H^{\prime}_{m+2}\}$ if $m>n-2$.
\begin{Lem}\label{lem-decomp-LA}
$L(\mathscr{G})_1=L(\mathscr{A})_1\sqcup\{H\cap H^{\prime}\mid H\in\mathscr{G},H^{\prime}\in\mathscr{G}\setminus\mathscr{A},H\neq H^{\prime}\}$.
Moreover if $H\in\mathscr{G}$ and $H^{\prime}\in\mathscr{G}\setminus\mathscr{A}$ with $H\neq H^{\prime}$, then $|\mathscr{G}_{H\cap H^{\prime}}|=2$.
\end{Lem}

\noindent
{\it Proof.}
Let $X\in L(\mathscr{G})_1$.
If $\mathscr{G}_X=\mathscr{A}_X$ then
$X=\bigcap_{H\in\mathscr{G}_X}H=\bigcap_{H\in\mathscr{A}_X}H
\in L(\mathscr{A})_1$.
Let $\mathscr{G}_X\neq\mathscr{A}_X$.
We suppose that $|\mathscr{G}_X|\geq 3$.
Then there exist $i\geq n+1$ and $X^{\prime}\in
L(\mathscr{A}\cup\{H^{\prime}_{n+1},\dots,H^{\prime}_{i-1}\})_1$
such that $X=X^{\prime}\cap H^{\prime}_i$.
By the condition (A), we have $X=\{0\}$.
This is a contradiction.
Therefore $|\mathscr{G}_X|=2$ and there exist
$H,H^{\prime}\in\mathscr{G}$
such that $X=H\cap H^{\prime}$.
Here if $H,H^{\prime}\in\mathscr{A}$ then $\mathscr{G}_X=\mathscr{A}_X$.
So either $H$ or $H^{\prime}$ does not belong to $\mathscr{A}$.
\hfill$\Box$
\vspace{2mm}

In the rest of this section, for a multi-set $\Psi$ and integers $e$ and $i$, $e^i\in\Psi$ means that the integer $e$ occurs $i$ times in $\Psi$.
\begin{Them}\label{thm-determine-exp}
Let $\ell=3$, let $m\geq n-2$, and
let $(\mathscr{A},V)$ be an essential arrangement.
Then $(\mathscr{A},V)$ is $m$-free with
\begin{align*}
{\rm exp}_{m}(\mathscr{A},V)=&
\{j+n-|\mathscr{A}_X|\mid X\in L(\mathscr{A})_1,
0\leq j\leq |\mathscr{A}_X|-2\}\\
&\cup\{(n-1)^{(m+2)n-\binom{n+1}{2}-
\sum_{X\in L(\mathscr{A})_1}(|\mathscr{A}_X|-1)},n^{\binom{m+2-n}{2}}\}.
\end{align*}
\end{Them}

\noindent
{\it Proof.}
By Theorem \ref{thm-D^m(A)=oplus_XD^i_Xdelta_X^m-i_X},
$(\mathscr{A},V)$ is $m$-free when $\ell=3$ and $m\geq n-2$.
By Lemma \ref{lem-decomp-LA},
\begin{align*}
L(\mathscr{G})_1=L(\mathscr{A})_1
\sqcup&\,\{H\cap H^{\prime}\mid
H,H^{\prime}\in\mathscr{G}\setminus\mathscr{A},H\neq H^{\prime}\}\\
\sqcup&\,\{H\cap H^{\prime}\mid
H\in\mathscr{A},H^{\prime}\in\mathscr{G}\setminus\mathscr{A}\}.
\end{align*}
Let $X\in L(\mathscr{A})_1$.
If $\mathscr{G}_X\neq\mathscr{A}_X$ then
$X\not\in L(\mathscr{A})_1$.
Hence $\mathscr{G}_X=\mathscr{A}_X$ and
$i_X=|\mathscr{G}_X|-2=|\mathscr{A}_X|-2\leq|\mathscr{A}_X|-1$.
By Proposition \ref{prop-2arrangements}, the multi-set of degrees of
the free $S$-basis for
$\bigoplus_{X\in L(\mathscr{A})_1}\bigoplus_{j=0}^{i_X}
SP_XD^{(j)}(\mathscr{A}_X,V_X)\delta_{X}^{m-j}$ is
\begin{align*}
&\bigcup_{X\in L(\mathscr{A})_1}\bigcup_{j=0}^{i_X}
\{j+n-|\mathscr{A}_X|,(n-|\mathscr{A}_X|+|\mathscr{A}_X|-1)^j\}\\
=&\bigcup_{X\in L(\mathscr{A})_1}
\left(\{j+n-|\mathscr{A}_X|\mid 0\leq j\leq |\mathscr{A}_X|-2\}
\cup\{(n-1)^{\sum_{j=0}^{|\mathscr{A}_X|-2}j}\}\right)\\
=&\{j+n-|\mathscr{A}_X|\mid
X\in L(\mathscr{A})_1,0\leq j\leq |\mathscr{A}_X|-2\}
\cup\{(n-1)^{\sum_{X\in L(\mathscr{A})_1}\binom{|\mathscr{A}_X|-1}{2}}\}\\
=&\{j+n-|\mathscr{A}_X|\mid
X\in L(\mathscr{A})_1,0\leq j\leq |\mathscr{A}_X|-2\}
\cup\{(n-1)^{\sum_{X\in L(\mathscr{A})_1}
\left(\binom{|\mathscr{A}_X|}{2}-(|\mathscr{A}_X|-1)\right)}\}\\
=&\{j+n-|\mathscr{A}_X|\mid
X\in L(\mathscr{A})_1,0\leq j\leq |\mathscr{A}_X|-2\}
\cup\{(n-1)^{\binom{n}{2}-\sum_{X\in L(\mathscr{A})_1}(|\mathscr{A}_X|-1)}\}.
\end{align*}
The last equation follows from Proposition \ref{prop-s_m=sum_Xs_iX} (2).

Let $\{H,H^{\prime}\}\subseteq\mathscr{G}\setminus\mathscr{A}$ and $X:=H\cap H^{\prime}$.
By the condition (A), we have $\mathscr{G}_X=\{H,H^{\prime}\}$, $i_X=0$, and $P_X=Q$. Therefore
\begin{align*}
\bigoplus_{X\in\{H\cap H^{\prime}\mid 
H,H^{\prime}\in\mathscr{G}\setminus\mathscr{A},H\neq H^{\prime}\}}
\bigoplus_{j=0}^{i_X}SP_XD^{(j)}(\mathscr{A}_X,V_X)\delta_{X}^{m-j}
=\bigoplus_{\{H,H^{\prime}\}\subseteq\mathscr{G}\setminus\mathscr{A}}
SQ\delta_{H\cap H^{\prime}}^{m}.
\end{align*}
The multi-set of degrees of the free $S$-basis for $\bigoplus_{\{H,H^{\prime}\}\subseteq\mathscr{G}\setminus\mathscr{A}}SQ\delta_{H\cap H^{\prime}}^{m}$ is $\{n^{\binom{\widetilde{n}-n}{2}}\}=\{n^{\binom{m+2-n}{2}}\}$.

Let $H\in\mathscr{A},H^{\prime}\in\mathscr{G}\setminus\mathscr{A}$
and $X:=H\cap H^{\prime}$.
Then by the condition (A), we have $\mathscr{G}_X=\{H,H^{\prime}\}$, $i_X=0$, and $P_X=Q/\alpha_H$. Therefore
\begin{align*}
\bigoplus_{X\in\{H\cap H^{\prime}\mid 
H\in\mathscr{A},H^{\prime}\in\mathscr{G}\setminus\mathscr{A}\}}
\bigoplus_{j=0}^{i_X}P_XSD^{(j)}(\mathscr{A}_X,V_X)\delta_{X}^{m-j}
=\bigoplus_{H\in\mathscr{A},H^{\prime}\in\mathscr{G}\setminus\mathscr{A}}
S(Q/\alpha_H)\delta_{H\cap H^{\prime}}^{m}.
\end{align*}
The multi-set of degrees of the free $S$-basis for $\bigoplus_{H\in\mathscr{A},H^{\prime}\in\mathscr{G}\setminus\mathscr{A}}S(Q/\alpha_H)\delta_{H\cap H^{\prime}}^{m}$ is $\{(n-1)^{(\widetilde{n}-n)n}\}=\{(n-1)^{(m+2)n-n^2}\}$.

By Theorem \ref{thm-D^m(A)=oplus_XD^i_Xdelta_X^m-i_X},
\begin{align*}
D^{(m)}(\mathscr{A},V)=&
\bigoplus_{X\in L(\mathscr{A})_1}\bigoplus_{j=0}^{i_X}
SP_X D^{(j)}(\mathscr{A}_X,V_X)\delta_X^{m-j}\\
&\oplus\bigoplus_{\{H,H^{\prime}\}\subseteq\mathscr{G}\setminus\mathscr{A}}
SQ\delta_{H\cap H^{\prime}}^{m}
\oplus\bigoplus_{H\in\mathscr{A},H^{\prime}\in\mathscr{G}\setminus\mathscr{A}}
S(Q/\alpha_H)\delta_{H\cap H^{\prime}}^{m}.
\end{align*}
Therefore
\begin{align*}
{\rm exp}_{m}(\mathscr{A},V)=&
\{j+n-|\mathscr{A}_X|\mid X\in L(\mathscr{A})_1,
0\leq j\leq |\mathscr{A}_X|-2\}\\
&\cup\{(n-1)^{(m+2)n-\binom{n+1}{2}-
\sum_{X\in L(\mathscr{A})_1}(|\mathscr{A}_X|-1)},n^{\binom{m+2-n}{2}}\}
\end{align*}
as required.
\hfill$\Box$

\section{Weyl arrangements of types A and B}\label{sec-wyle-arrAB}
Weyl arrangements $(\mathcal{A}_{\ell},V)$ and $(\mathcal{B}_{\ell},V)$ are defined by
\begin{align*}
\mathcal{A}_{\ell}&:=\left\{\{x_i=0\}\,\middle|\,i=1,\dots,\ell\right\}\cup\left\{\{x_i-x_j=0\}\,\middle|\,1\leq i<j\leq \ell\right\},\\
\mathcal{B}_{\ell}&:=\left\{\{x_i=0\}\,\middle|\,i=1,\dots,\ell\right\}\cup\left\{\{x_i-x_j=0\},\{x_i+x_j=0\}\,\middle|\,1\leq i<j\leq \ell\right\}.
\end{align*}
The author \cite{N-Cox2free} proved that $\mathcal{A}_{\ell}$ and $\mathcal{B}_{\ell}$ are $2$-free for all $\ell\geq 1$.
By Theorem \ref{thm-D^m(A)=oplus_XD^i_Xdelta_X^m-i_X}, when $\ell=3$, $\mathcal{A}_{3}$ is $m$-free if $m\geq|\mathcal{A}_{3}|-2=4$, while $\mathcal{B}_{3}$ is $m^{\prime}$-free if $m^{\prime}\geq|\mathcal{B}_{3}|-2=7$.
Therefore, to prove that $\mathcal{A}_{3}$ and $\mathcal{B}_{3}$ are $m$-free for all $m\geq 0$, it is enough to prove that $\mathcal{A}_{3}$ is $3$-free and $\mathcal{B}_{3}$ is $i$-free for any $i\in\{3,4,5,6\}$.

We define an operator
\begin{align}\label{eq-def-thetak-weyl}
\theta_{k}:=\sum_{|\bm{a}|=m}\left(\sum_{r=1}^{\ell}a_rx^k_r\right)\frac{x^{\bm{a}}}{\bm{a}!}\partial^{\bm{a}}
\end{align}
for $k\geq 0$, where we write $\bm{a}=(a_1,\dots,a_{\ell})$ for any $|\bm{a}|=m$ in \eqref{eq-def-thetak-weyl}.
We note that $\deg(\theta_k)=m+k$.
\begin{Lem}\label{lem-thetak-in-DADB}
Let $\ell\geq 1$. Then the following hold.
\begin{itemize}
\item[(1)] $\theta_k\in D^{(m)}(\mathcal{A}_{\ell})$ for any $k\geq 0$.
\item[(2)] If $k$ is even, then $\theta_k\in D^{(m)}(\mathcal{B}_{\ell})$.
\end{itemize}
\end{Lem}

\noindent
{\it Proof.}
(1)\ Let $\bm{b}\in\mathbb{N}^{\ell}$ with $|\bm{b}|=m-1$.
For any $1\leq i\leq\ell$, we have
\begin{align*}
\theta_{k}\left(x_ix^{\bm{b}}\right)=x_ix^{\bm{b}}\left(\sum_{r=1}^{\ell}b_r x_r^k+x_i^k\right)\in x_i S.
\end{align*}
Meanwhile for any $i,j$ with $1\leq i<j\leq \ell$,
\begin{align*}
\theta_{k}\left((x_i-x_j)x^{\bm{b}}\right)
&=x_ix^{\bm{b}}\left(\sum_{r=1}^{\ell}b_r x_r^k+x_i^k\right)
-x_jx^{\bm{b}}\left(\sum_{r=1}^{\ell}b_r x_r^k+x_j^k\right)\\
&=x^{\bm{b}}\left((x_i-x_j)\sum_{r=1}^{\ell}b_r x_r^k+x_i^{k+1}-x_j^{k+1}\right)\\
&=(x_i-x_j)x^{\bm{b}}\left(\sum_{r=1}^{\ell}b_r x_r^k+\sum_{t=0}^{k}x_i^{k-t} x_j^{t}\right)\in (x_i-x_j)S.
\end{align*}
Therefore $\theta_k\in D^{(m)}(\mathcal{A}_{\ell},V)$ by Proposition \ref{prop-check-opes}.

(2)\ Let $\bm{b}\in\mathbb{N}^{\ell}$ with $|\bm{b}|=m-1$.
Let $i,j$ be integers with $1\leq i<j\leq \ell$.
Since $k$ is even, we have $x_{i}^{k+1}+x_{j}^{k+1}=(x_i+x_j)\sum_{t=0}^{k}(-1)^t x_i^{k-t} x_j^{t}$.
Then
\begin{align*}
\theta_{k}\left((x_i+x_j)x^{\bm{b}}\right)
&=x_ix^{\bm{b}}\left(\sum_{r=1}^{\ell}b_r x_r^k+x_i^k\right)
+x_jx^{\bm{b}}\left(\sum_{r=1}^{\ell}b_r x_r^k+x_j^k\right)\\
&=x^{\bm{b}}\left((x_i+x_j)\sum_{r=1}^{\ell}b_r x_r^k+x_i^{k+1}+x_j^{k+1}\right)\\
&=(x_i+x_j)x^{\bm{b}}\left(\sum_{r=1}^{\ell}b_r x_r^k+\sum_{t=0}^{k}(-1)^t x_i^{k-t} x_j^{t}\right)\in (x_i+x_j)S.
\end{align*}
Combining the arguments of (1), we have $\theta_k\in D^{(m)}(\mathcal{B}_{\ell},V)$ by Proposition \ref{prop-check-opes}.
\hfill$\Box$

\subsection{The arrangement $\mathcal{A}_3$}\label{subsec-BasisA3}
Let $\ell=3$ and we consider the arrangement $\mathscr{A}=\mathcal{A}_3$. The operator $\theta_k$ belongs to $D^{(3)}(\mathcal{A}_3,V)$ by Lemma \ref{lem-thetak-in-DADB}, while we can prove similarly to the proof of Theorem \ref{thm-D^m(A)=oplus_XD^i_Xdelta_X^m-i_X} that $P_X\theta\delta_X^{3-j}\in D^{(3)}(\mathcal{A}_3,V)$ for any $X\in L(\mathcal{A}_3)_1$, for any $0\leq j\leq i_X$, and for any $\theta\in D^{(j)}\left((\mathcal{A}_3)_X,V_X\right)$.
We construct a free basis for $D^{(3)}(\mathcal{A}_3,V)$ consisting of these operators.
We have exactly seven elements in $L(\mathcal{A}_3)_1$:
\begin{align*}
X_{1}&=\{x_2=0\}\cap\{x_3=0\}\cap\{x_2-x_3=0\},\\
X_{2}&=\{x_1=0\}\cap\{x_3=0\}\cap\{x_1-x_3=0\},\\
X_{3}&=\{x_1=0\}\cap\{x_2=0\}\cap\{x_1-x_2=0\},\\
X_{4}&=\{x_1-x_2=0\}\cap\{x_1-x_3=0\}\cap\{x_2-x_3=0\},\\
X_{5}&=\{x_3=0\}\cap\{x_1-x_2=0\},\\
X_{6}&=\{x_2=0\}\cap\{x_1-x_3=0\},\\
X_{7}&=\{x_1=0\}\cap\{x_2-x_3=0\}.
\end{align*}
Then
\begin{align*}
&P_{X_1}=x_1(x_1-x_2)(x_1-x_3),\ P_{X_2}=x_2(x_1-x_2)(x_2-x_3),\ P_{X_3}=x_3(x_1-x_3)(x_2-x_3),\\
&P_{X_4}=x_1x_2x_3,\ P_{X_5}=x_1x_2(x_1-x_3)(x_2-x_3),\ P_{X_6}=x_1x_3(x_1-x_2)(x_2-x_3),\\
&P_{X_7}=x_2x_3(x_1-x_2)(x_1-x_3),
\end{align*}
while
\begin{align*}
&\delta_{X_1}=\partial_1,\ \delta_{X_2}=\partial_2,\ \delta_{X_3}=\partial_3,\ \delta_{X_4}=\partial_1+\partial_2+\partial_3,\\
&\delta_{X_5}=\partial_1+\partial_2,\ \delta_{X_6}=\partial_1+\partial_3,\ \delta_{X_7}=\partial_2+\partial_3.
\end{align*}
\begin{Prop}\label{prop-BasisA3-m=3}
The set
\begin{align}\label{eq-basis-A3m3}
\left\{P_{X_j}\delta_{X_j}^3\,\middle|\,1\leq j\leq 7\right\}\cup\left\{\theta_0, \theta_1, \theta_2\right\}
\end{align}
is a free basis for $D^{(3)}(\mathcal{A}_3,V)$, while ${\rm exp}_3(\mathcal{A}_3,V)=\{3^5, 4^4, 5^1\}$, where $e^i$ means that the integer $e$ occurs $i$ times in the multi-set.
\end{Prop}

\noindent
{\it Proof.}
We prove the assertion by Saito's criterion.
We have that $s_2(3)=\binom{2+2}{2}=6$ and $\prod_{i=1}^{7}P_{X_i}=x_1^4 x_2^4 x_3^4(x_1-x_2)^4(x_1-x_3)^4(x_2-x_3)^4=Q(\mathcal{A}_3)^4$.
It is enough to prove that $\det M_{3}(\delta_{X_1}^3, \delta_{X_2}^3, \delta_{X_3}^3, \delta_{X_4}^3, \delta_{X_5}^3, \delta_{X_6}^3, \delta_{X_7}^3, \theta_0, \theta_1, \theta_2)=cQ(\mathcal{A}_3)^2$ for some $c\in\mathbb{K}\setminus\{0\}$.
Since $\delta_{X_1}^3=\partial_1^3$, $\delta_{X_2}^3=\partial_2^3$, and $\delta_{X_3}^3=\partial_3^3$, we can omit coefficients of $\partial_1^3,\partial_2^3,\partial_3^3$ in the next determinant calculations.
We use $\doteq$ to denote equality up to multiplication by a nonzero element of $\mathbb{K}$.
Then
\begin{align}
&\det\left(M_{3}(\delta_{X_1}^3, \delta_{X_2}^3, \delta_{X_3}^3, \delta_{X_4}^3, \delta_{X_5}^3, \delta_{X_6}^3, \delta_{X_7}^3, \theta_0, \theta_1, \theta_2)\right)\notag\\
\doteq&\,\left|
\begin{matrix}
3&1&0&0&x_1^2x_2&x_1^2x_2(2x_1+x_2)&x_1^2x_2(2x_1^2+x_2^2)\\
3&0&1&0&x_1^2x_3&x_1^2x_3(2x_1+x_3)&x_1^2x_3(2x_1^2+x_3^2)\\
3&1&0&0&x_2^2x_1&x_2^2x_1(2x_2+x_1)&x_2^2x_1(2x_2^2+x_1^2)\\
3&0&0&1&x_2^2x_3&x_2^2x_3(2x_2+x_3)&x_2^2x_3(2x_2^2+x_3^2)\\
3&0&1&0&x_3^2x_1&x_3^2x_1(2x_3+x_1)&x_3^2x_1(2x_3^2+x_1^2)\\
3&0&0&1&x_3^2x_2&x_3^2x_2(2x_3+x_2)&x_3^2x_2(2x_3^2+x_2^2)\\
6&0&0&0&x_1 x_2 x_3&x_1 x_2 x_3(x_1+x_2+x_3)&x_1 x_2 x_3(x_1^2+x_2^2+x_3^2)
\end{matrix}
\right|\tag{$\ast$}\label{eq-det-coef-D3A3}\\
\doteq&\,\left|
\begin{matrix}
1&0&0&x_1^2x_2&x_1^2x_2(2x_1+x_2)&x_1^2x_2(2x_1^2+x_2^2)\\
0&1&0&x_1^2x_3&x_1^2x_3(2x_1+x_3)&x_1^2x_3(2x_1^2+x_3^2)\\
1&0&0&x_2^2x_1&x_2^2x_1(2x_2+x_1)&x_2^2x_1(2x_2^2+x_1^2)\\
0&0&1&x_2^2x_3&x_2^2x_3(2x_2+x_3)&x_2^2x_3(2x_2^2+x_3^2)\\
0&1&0&x_3^2x_1&x_3^2x_1(2x_3+x_1)&x_3^2x_1(2x_3^2+x_1^2)\\
0&0&1&x_3^2x_2&x_3^2x_2(2x_3+x_2)&x_3^2x_2(2x_3^2+x_2^2)
\end{matrix}
\right|\notag\\
\doteq&\,\left|
\begin{matrix}
x_2x_1(x_2-x_1)&x_2x_1(2x_2^2-2x_1^2)&x_2x_1(2x_2^3+x_2x_1^2-x_2^2x_1-2x_1^3)\\
x_3x_1(x_3-x_1)&x_3x_1(2x_3^2-2x_1^2)&x_3x_1(2x_3^3+x_3x_1^2-x_3^2x_1-2x_1^3)\\
x_3x_2(x_3-x_2)&x_3x_2(2x_3^2-2x_2^2)&x_3x_2(2x_3^3+x_3x_2^2-x_3^2x_2-2x_2^3)
\end{matrix}
\right|\notag\\
=&\,x_1^2x_2^2x_3^2(x_2-x_1)(x_3-x_1)(x_3-x_2)\left|
\begin{matrix}
1&2(x_2+x_1)&2x_2^2+x_2x_1+2x_1^2\\
1&2(x_3+x_1)&2x_3^2+x_3x_1+2x_1^2\\
1&2(x_3+x_2)&2x_3^2+x_3x_2+2x_2^2
\end{matrix}
\right|\notag\\
\doteq&\,x_1^2x_2^2x_3^2(x_2-x_1)(x_3-x_1)(x_3-x_2)\left|
\begin{matrix}
x_3-x_2&2x_3^2-2x_2^2+x_1(x_3-x_2)\\
x_3-x_1&2x_3^2-2x_1^2+x_2(x_3-x_1)
\end{matrix}
\right|\notag\\
=&\,x_1^2x_2^2x_3^2(x_2-x_1)(x_3-x_1)^2(x_3-x_2)^2\left|
\begin{matrix}
1&2x_3+2x_2+x_1\\
1&2x_3+2x_1+x_2
\end{matrix}
\right|\doteq Q(\mathcal{A}_3)^2,\notag
\end{align}
where rows are arranged in the coefficients of $\partial_1^2\partial_2, \partial_1^2\partial_3, \partial_2^2\partial_1, \partial_2^2\partial_3, \partial_3^2\partial_1, \partial_3^2\partial_2, \partial_1\partial_2\partial_3$ in the determinant \eqref{eq-det-coef-D3A3}.
Therefore the set \eqref{eq-basis-A3m3} forms a basis for $D^{(3)}(\mathcal{A}_3,V)$.
\hfill$\Box$
\begin{Cor}
The arrangement $\mathcal{A}_3$ is $m$-free for any $m\geq 0$.
\end{Cor}

\subsection{The arrangement $\mathcal{B}_3$}\label{subsec-BasisB3}
Let $\ell=3$ and we consider the arrangement $\mathscr{A}=\mathcal{B}_3$. If $k$ is even, then the operator $\theta_k$ belongs to $D^{(m)}(\mathcal{B}_3,V)$ by Lemma \ref{lem-thetak-in-DADB}. Meanwhile $P_X\theta\delta_X^{m-j}\in D^{(m)}(\mathcal{B}_3,V)$ for any $X\in L(\mathcal{B}_3)_1$, for any $0\leq j\leq i_X$, and for any $\theta\in D^{(j)}\left((\mathcal{B}_3)_X,V_X\right)$.
We have exactly thirteen elements in $L(\mathcal{B}_3)_1$:
\begin{align*}
X_{1}&=\{x_2=0\}\cap\{x_3=0\}\cap\{x_2-x_3=0\}\cap\{x_2+x_3=0\},\\
X_{2}&=\{x_1=0\}\cap\{x_3=0\}\cap\{x_1-x_3=0\}\cap\{x_1+x_3=0\},\\
X_{3}&=\{x_1=0\}\cap\{x_2=0\}\cap\{x_1-x_2=0\}\cap\{x_1+x_2=0\},\\
X_{4}&=\{x_1-x_2=0\}\cap\{x_1-x_3=0\}\cap\{x_2-x_3=0\},\\
X_{5}&=\{x_1+x_2=0\}\cap\{x_1+x_3=0\}\cap\{x_2-x_3=0\},\\
X_{6}&=\{x_1+x_2=0\}\cap\{x_1-x_3=0\}\cap\{x_2+x_3=0\},\\
X_{7}&=\{x_1-x_2=0\}\cap\{x_1+x_3=0\}\cap\{x_2+x_3=0\},\\
X_{8}&=\{x_3=0\}\cap\{x_1-x_2=0\},\ 
X_{9}=\{x_2=0\}\cap\{x_1-x_3=0\},\\
X_{10}&=\{x_1=0\}\cap\{x_2-x_3=0\},\ 
X_{11}=\{x_3=0\}\cap\{x_1+x_2=0\},\\
X_{12}&=\{x_2=0\}\cap\{x_1+x_3=0\},\ 
X_{13}=\{x_1=0\}\cap\{x_2+x_3=0\}.
\end{align*}
Then
\begin{align*}
&P_{X_1}=x_1(x_1^2-x_2^2)(x_1^2-x_3^2),\ P_{X_2}=x_2(x_1^2-x_2^2)(x_2^2-x_3^2),\ P_{X_3}=x_3(x_1^2-x_3^2)(x_2^2-x_3^2),\\
&P_{X_4}=x_1x_2x_3(x_1+x_2)(x_1+x_3)(x_2+x_3),\ 
P_{X_5}=x_1x_2x_3(x_1-x_2)(x_1-x_3)(x_2+x_3),\\
&P_{X_6}=x_1x_2x_3(x_1-x_2)(x_1+x_3)(x_2-x_3),\ 
P_{X_7}=x_1x_2x_3(x_1+x_2)(x_1-x_3)(x_2-x_3),\\
&P_{X_8}=x_1x_2(x_1+x_2)(x_1^2-x_3^2)(x_2^2-x_3^2),\ 
P_{X_9}=x_1x_3(x_1^2-x_2^2)(x_1+x_3)(x_2^2-x_3^2),\\
&P_{X_{10}}=x_2x_3(x_1^2-x_2^2)(x_1^2-x_3^2)(x_2+x_3),\ 
P_{X_{11}}=x_1x_2(x_1-x_2)(x_1^2-x_3^2)(x_2^2-x_3^2),\\
&P_{X_{12}}=x_1x_3(x_1^2-x_2^2)(x_1-x_3)(x_2^2-x_3^2),\ 
P_{X_{13}}=x_2x_3(x_1^2-x_2^2)(x_1^2-x_3^2)(x_2-x_3),
\end{align*}
while
\begin{align*}
&\delta_{X_1}=\partial_1, \delta_{X_2}=\partial_2, \delta_{X_3}=\partial_3,\ 
\delta_{X_4}=\partial_1+\partial_2+\partial_3,\ 
\delta_{X_5}=-\partial_1+\partial_2+\partial_3,\\
&\delta_{X_6}=\partial_1-\partial_2+\partial_3,\ 
\delta_{X_7}=\partial_1+\partial_2-\partial_3,\ 
\delta_{X_8}=\partial_1+\partial_2,\ 
\delta_{X_9}=\partial_1+\partial_3,\\
&\delta_{X_{10}}=\partial_2+\partial_3,\ 
\delta_{X_{11}}=\partial_1-\partial_2,\ 
\delta_{X_{12}}=\partial_1-\partial_3,\ 
\delta_{X_{13}}=\partial_2-\partial_3.
\end{align*}
In addition we define Euler derivations
\begin{align*}
&\varepsilon_{X_1}=x_2\partial_2+x_3\partial_3,\ 
\varepsilon_{X_2}=x_1\partial_1+x_3\partial_3,\ 
\varepsilon_{X_3}=x_1\partial_1+x_2\partial_2, \\
&\varepsilon_{X_4}=\frac{x_1-x_2}{3}\left(\partial_1-2\partial_2+\partial_3\right)+\frac{x_1-x_3}{3}\left(\partial_1+\partial_2-2\partial_3\right),\\
&\varepsilon_{X_5}=\frac{x_1+x_3}{3}\left(\partial_1-\partial_2+2\partial_3\right)+\frac{x_1+x_2}{3}\left(\partial_1+2\partial_2-\partial_3\right),\\
&\varepsilon_{X_6}=\frac{x_1+x_2}{3}\left(2\partial_1+\partial_2-\partial_3\right)+\frac{x_2+x_3}{3}\left(-\partial_1+\partial_2+2\partial_3\right),\\
&\varepsilon_{X_7}=\frac{x_2+x_3}{3}\left(-\partial_1+2\partial_2+\partial_3\right)+\frac{x_1+x_3}{3}\left(2\partial_1-\partial_2+\partial_3\right),\\
&\varepsilon_{X_8}=\frac{x_1-x_2}{2}\left(\partial_1-\partial_2\right)+x_3\partial_3
\end{align*}
with respect to $X_j\ (1\leq j\leq 8)$.
Then $\varepsilon_{X_j}\in D^{(1)}\left(\left(\mathcal{B}_3\right)_{X_{j}},V_{X_j}\right)$ for any $1\leq j\leq 8$.
Although there are thirteen Euler derivations with respect to $X\in L(\mathcal{B}_3)_1$, we only use these eight operators to construct free bases for $D^{(i)}(\mathcal{B}_3,V)$ for $i\in\{3,4,5,6\}$.
\begin{Prop}\label{prop-basis-DB3m=3-6}
In the following multi-sets of exponents, $e^i$ means that the integer $e$ occurs $i$ times.

\noindent
(1)\ The set
\begin{align}\label{eq-basis-B3m3}
\left\{P_{X_j}\delta_{X_j}^3\,\middle|\, 1\leq j\leq 4\right\}\cup\left\{P_{X_j}\varepsilon_{X_j}\delta_{X_j}^2\,\middle|\, 1\leq j\leq 4\right\}\cup\left\{\theta_0, \theta_2\right\}
%\{\theta_0, \theta_2, P_{X_1}\delta_{X_1}^3, P_{X_2}\delta_{X_2}^3, P_{X_3}\delta_{X_3}^3,P_{X_4}\delta_{X_4}^3, P_{X_1}\varepsilon_{X_1}\delta_{X_1}^2, P_{X_2}\varepsilon_{X_2}\delta_{X_2}^2, P_{X_3}\varepsilon_{X_3}\delta_{X_3}^2, P_{X_4}\varepsilon_{X_4}\delta_{X_4}^2\}
\end{align}
is a free basis for $D^{(3)}(\mathcal{B}_3,V)$, while ${\rm exp}_3(\mathcal{B}_3,V)=\{3^1, 5^4, 6^4, 7^1\}$.

\noindent
(2)\ The set
\begin{align}\label{eq-basis-B3m4}
\left\{P_{X_j}\delta_{X_j}^4\,\middle|\, 1\leq j\leq 10\right\}\cup\left\{P_{X_1}\varepsilon_{X_1}\delta_{X_1}^3, P_{X_2}\varepsilon_{X_2}\delta_{X_2}^3, P_{X_3}\varepsilon_{X_3}\delta_{X_3}^3, P_{X_8}\varepsilon_{X_8}\delta_{X_8}^3\right\}\cup\left\{\theta_0\right\}
%\left\{
%\begin{array}{l}
%\theta_0, P_{X_1}\delta_{X_1}^4, P_{X_2}\delta_{X_2}^4, P_{X_3}\delta_{X_3}^4, P_{X_4}\delta_{X_4}^4, P_{X_5}\delta_{X_5}^4, P_{X_6}\delta_{X_6}^4, P_{X_7}\delta_{X_7}^4, P_{X_8}\delta_{X_8}^4,\\
%P_{X_9}\delta_{X_9}^4, P_{X_{10}}\delta_{X_{10}}^4, P_{X_1}\varepsilon_{X_1}\delta_{X_1}^3, P_{X_2}\varepsilon_{X_2}\delta_{X_2}^3, P_{X_3}\varepsilon_{X_3}\delta_{X_3}^3, P_{X_8}\varepsilon_{X_8}\delta_{X_8}^3
%\end{array}
%\right\}
\end{align}
is a free basis for $D^{(4)}(\mathcal{B}_3,V)$, while ${\rm exp}_4(\mathcal{B}_3,V)=\{4^1, 5^3, 6^7, 7^3, 8^1\}$.

\noindent
(3)\ The set
\begin{align}\label{eq-basis-B3m5}
\begin{array}{l}
\left\{P_{X_j}\delta_{X_j}^5\,\middle|\, 1\leq j\leq 13\right\}\cup\left\{P_{X_j}\varepsilon_{X_j}\delta_{X_j}^4\,\middle|\, 1\leq j\leq 4\right\}\cup\left\{\theta_0\right\}\\
\cup\left\{P_{X_1}x_2(x_2^2-x_3^2)\partial_2\delta_{X_1}^4, P_{X_2}x_3(x_1^2-x_3^2)\partial_3\delta_{X_2}^4, P_{X_3}x_1(x_1^2-x_2^2)\partial_1\delta_{X_3}^4\right\}
\end{array}
%\left\{
%\begin{array}{l}
%\theta_0, P_{X_1}\delta_{X_1}^5, P_{X_2}\delta_{X_2}^5, P_{X_3}\delta_{X_3}^5, P_{X_4}\delta_{X_4}^5, P_{X_5}\delta_{X_5}^5, P_{X_6}\delta_{X_6}^5, P_{X_7}\delta_{X_7}^5, P_{X_8}\delta_{X_8}^5, P_{X_9}\delta_{X_9}^5,\\
%P_{X_{10}}\delta_{X_{10}}^5, P_{X_{11}}\delta_{X_{11}}^5, P_{X_{12}}\delta_{X_{12}}^5, P_{X_{13}}\delta_{X_{13}}^5, P_{X_1}\varepsilon_{X_1}\delta_{X_1}^4, P_{X_2}\varepsilon_{X_2}\delta_{X_2}^4, P_{X_3}\varepsilon_{X_3}\delta_{X_3}^4,\\
%P_{X_4}\varepsilon_{X_4}\delta_{X_4}^4, P_{X_1}x_2(x_2^2-x_3^2)\partial_2\delta_{X_1}^4, P_{X_2}x_3(x_1^2-x_3^2)\partial_3\delta_{X_2}^4, P_{X_3}x_1(x_1^2-x_2^2)\partial_1\delta_{X_3}^4
%\end{array}
%\right\}
\end{align}
is a free basis for $D^{(5)}(\mathcal{B}_3,V)$, while ${\rm exp}_5(\mathcal{B}_3,V)=\{5^4, 6^7, 7^7, 8^3\}$.

\noindent
(4)\ The set
\begin{align}\label{eq-basis-B3m6}
\begin{array}{l}
\left\{P_{X_j}\delta_{X_j}^6\,\middle|\, 1\leq j\leq 13\right\}\cup\left\{P_{X_j}\varepsilon_{X_j}\delta_{X_j}^5\,\middle|\, 1\leq j\leq 7\right\}\cup\left\{\theta_0\right\}\\
\cup\left\{
\begin{array}{l}
P_{X_1}x_2(x_2^2-x_3^2)\partial_2\delta_{X_1}^5,\\
P_{X_2}x_3(x_1^2-x_3^2)\partial_3\delta_{X_2}^5,\\
P_{X_3}x_1(x_1^2-x_2^2)\partial_1\delta_{X_3}^5,\\
P_{X_4}(x_1-x_3)(x_2-x_3)(\partial_1+\partial_2-2\partial_3)\delta_{X_4}^5,\\
P_{X_5}(x_1+x_2)(x_2-x_3)(\partial_1+2\partial_2-\partial_3)\delta_{X_5}^5,\\ 
P_{X_6}(x_1-x_3)(x_2+x_3)(-\partial_1+\partial_2+2\partial_3)\delta_{X_6}^5,\\
P_{X_7}(x_1-x_2)(x_1+x_3)(2\partial_1-\partial_2+\partial_3)\delta_{X_7}^5
\end{array}
\right\}
\end{array}
%\left\{
%\begin{array}{ll}
%\theta_0, 
%P_{X_1}\delta_{X_1}^6, 
%P_{X_2}\delta_{X_2}^6, 
%P_{X_3}\delta_{X_3}^6, 
%P_{X_4}\delta_{X_4}^6, 
%P_{X_5}\delta_{X_5}^6, 
%P_{X_6}\delta_{X_6}^6, 
%P_{X_7}\delta_{X_7}^6, 
%P_{X_8}\delta_{X_8}^6, 
%P_{X_9}\delta_{X_9}^6,
%P_{X_{10}}\delta_{X_{10}}^6,\\
%P_{X_{11}}\delta_{X_{11}}^6, 
%P_{X_{12}}\delta_{X_{12}}^6, 
%P_{X_{13}}\delta_{X_{13}}^6, 
%P_{X_1}\varepsilon_{X_1}\delta_{X_1}^5, 
%P_{X_2}\varepsilon_{X_2}\delta_{X_2}^5, 
%P_{X_3}\varepsilon_{X_3}\delta_{X_3}^5, 
%P_{X_4}\varepsilon_{X_4}\delta_{X_4}^5, 
%P_{X_5}\varepsilon_{X_5}\delta_{X_5}^5,\\
%P_{X_6}\varepsilon_{X_6}\delta_{X_6}^5, 
%P_{X_7}\varepsilon_{X_7}\delta_{X_7}^5, 
%P_{X_1}x_2(x_2^2-x_3^2)\partial_2\delta_{X_1}^5, 
%P_{X_2}x_3(x_1^2-x_3^2)\partial_3\delta_{X_2}^5, 
%P_{X_3}x_1(x_1^2-x_2^2)\partial_1\delta_{X_3}^5,\\
%P_{X_4}(x_1-x_3)(x_2-x_3)(\partial_1+\partial_2-2\partial_3)\delta_{X_4}^5, 
%P_{X_5}(x_1+x_2)(x_2-x_3)(\partial_1+2\partial_2-\partial_3)\delta_{X_5}^5,\\ 
%P_{X_6}(x_1-x_3)(x_2+x_3)(-\partial_1+\partial_2+2\partial_3)\delta_{X_6}^5,
%P_{X_7}(x_1-x_2)(x_1+x_3)(2\partial_1-\partial_2+\partial_3)\delta_{X_7}^5
%\end{array}
%\right\}
\end{align}
is a free basis for $D^{(6)}(\mathcal{B}_3,V)$, while ${\rm exp}_6(\mathcal{B}_3,V)=\{5^3, 6^8, 7^{10}, 8^7\}$.
\end{Prop}
\vspace{2mm}

We can prove Proposition \ref{prop-basis-DB3m=3-6} by straightforward determinant calculations of coefficient matrices.
However it is too long to write them all down.
We omit the proofs of Proposition \ref{prop-basis-DB3m=3-6}.
\begin{Cor}
The arrangement $\mathcal{B}_3$ is $m$-free for any $m\geq 0$.
\end{Cor}

\section*{Acknowledgments}
This work was supported by JSPS KAKENHI Grant Number 20K14282.

\end{document}